\documentclass[draft,reqno,12pt]{amsart}

\newtheorem{theorem}{Theorem}[section]
\newtheorem{proposition}[theorem]{Proposition}
\newtheorem{lemma}[theorem]{Lemma}
\newtheorem{cor}[theorem]{Corollary}

\theoremstyle{definition}

\newtheorem*{rem}{[Remark]}
\newtheorem{exam}[theorem]{Example}

\newcommand{\Gin}{\ensuremath{\mathrm{Gin}}}
\newcommand{\gin}{\ensuremath{\mathrm{gin}}}

\newcommand{\Gins}{\ensuremath{\mathrm{Gins}}}
\newcommand{\gins}{\ensuremath{\mathrm{gins}}}

\newcommand{\init}{\ensuremath{\mathrm{in}}}
\newcommand{\GL}{\ensuremath{GL_n (K)}}

\newcommand{\reg}{\ensuremath{\mathrm{reg}}}

\newcommand{\dele}[1]{\ensuremath{\Delta^e({#1})}}

\newcommand{\lex}{{<_{\mathrm{lex}}}}
\newcommand{\rlex}{{<_{\mathrm{rev}}}}

\newcommand{\T}{ \mathrm{Trans}}
\newcommand{\F}{ \mathcal{F}}
\def\cocoa{{\hbox{\rm C\kern-.13em o\kern-.07em C\kern-.13em o\kern-.15em A}}}

\begin{document}

\title[Monomial ideals whose generic initial ideal is unique]
{
Monomial ideals arising from flag complexes whose generic initial ideals do not depend on term orders}
\author{Satoshi Murai}
\address{
Department of Pure and Applied Mathematics\\
Graduate School of Information Science and Technology\\
Osaka University\\
Toyonaka, Osaka, 560-0043, Japan\\
}
\email{s-murai@ist.osaka-u.ac.jp
}
\thanks{The author is supported by JSPS Research Fellowships for Young Scientists}

\begin{abstract}
We will study monomial ideals $I$ in the exterior algebra
as well as in the polynomial ring whose generic initial ideal
is constant for all term orders
up to permutations of variables.
First, in the exterior algebra,
we determine all graphs and all flag complexes whose 
exterior face ideal satisfies the above condition.
Second, in the polynomial ring, it will be shown that
the generic initial ideal $\gin_\sigma(I(G))$ of the edge ideal $I(G)$ of a graph $G$
is constant for all term orders $\sigma$ up to permutations of variables
if and only if $G$ is a complete bipartite graph.
\end{abstract}

\maketitle

\section{Introduction}
In this paper,
we will study some monomial ideals $J$ in the exterior algebra $E$
which satisfies $\Gin_\sigma(J)=\Gin_\tau(J)$ 
for all term orders $\sigma$ and $\tau$ satisfying $e_1 > e_2 > \cdots > e_n$,
where $\Gin_\sigma(I)$ is the generic initial ideal of $J$
with respect to the term order $\sigma$
and where $E=\bigoplus_{k=0}^n \bigwedge^k V$
is the exterior algebra of a vector space $V$ over an infinite field $K$
with basis $e_1,\dots,e_n$.
Our results determine all graphs and all flag complexes
whose exterior face ideal satisfies the above condition.

Let $\Gamma$ be a simplicial complex on $[n]=\{1,\dots,n\}$.
For a subset $S=\{i_1,\dots,i_k\} \subset [n]$,
the element $e_S=e_{i_1} \wedge \cdots \wedge e_{i_k} \in E$ will be called a \textit{monomial} in $E$ of degree $k$,
where $i_1 < \cdots < i_k$.
The \textit{exterior face ideal $J_\Gamma$ of} $\Gamma$ is the monomial ideal in $E$
generated by all monomials $e_S \in E$ with $S \not \in \Gamma$.
In this paper,
a $1$-dimensional simplicial complex will be called a \textit{graph}.
If $G$ is a graph on $[n]$,
then we write $\F(G)$ for the flag complex of $G$,
where the flag complex of $G$ is the simplicial complex defined by
$$ \F(G) = \{ S \subset [n]: \{i,j\} \in G\ \mbox{ for all } \{i,j\} \subset S\}.$$
For a graded ideal $J$ in $E$,
let
$$ \Gins (J) = \{ \Gin_\sigma(J): \sigma \mbox{ is a term order with } e_1 >_\sigma e_2>_\sigma \cdots >_\sigma e_n\}.$$
We will consider graphs $G$ which satisfy $|\Gins(J_{\F(G)})|=1$,
where $|A|$ denotes the cardinality of a finite set $A$.

A similar problem was asked by Kalai in \cite[Problem 7]{K}.
Let $J \subset E$ be a homogeneous ideal
and $\init_\sigma(J)$ the initial ideal of $J$ with respect to a term order $\sigma$.
A monomial ideal $J \subset E$ is said to be \textit{strongly stable}
if $e_S \in J$ implies $e_{(S \setminus \{j\}) \cup \{i\}} \in J$ for all $j \in S$ and $i \not \in S$ with $i<j$.
Kalai asked which simplicial complex $\Gamma$ satisfies
$\init_{\mathrm{rev}} (\varphi (J_\Gamma))= \gin(J_\Gamma)$ for all $\varphi \in \GL$ with which $\init_{\mathrm{rev}}( \varphi(J_\Gamma))$
is strongly stable,
where the general linear group $\GL$ acts linearly on $E$
and $\init_{\mathrm{rev}}(J)$ is the initial ideal of $J \subset E$ with respect to the reverse lexicographic order. 
We will show that
if $G$ satisfies $|\Gins(J_G)|=1$ (resp. $|\Gins(J_{\F(G)})|=1$) then $G$ (resp. $\F(G)$) satisfies the above condition.

A strongly stable ideal $J' \subset E$ is called a \textit{transformed strongly stable ideal of $J$}
if there exist $\varphi_1,\varphi_2,\dots,\varphi_k \in \GL$ and term orders
$\sigma_1,\sigma_2,\dots,\sigma_k$ such that
$$J'= \init_{\sigma_k} ( \varphi_k ( \cdots \init_{\sigma_{2}} ( \varphi_2 (\init_{\sigma_1}(\varphi_1 (J)) )) \cdots )).$$
Set
$$\T (J)= \{ J'\subset E: J' \mbox{ is a transformed strongly stable ideal of } J\}.$$
Since every generic initial ideal is strongly stable,
$\T(J)$ contains $\Gins(J)$ together with all strongly stable ideals considered in \cite[Problem 7]{K}.
However,
we will show that $|\T(J_{\F(G)})|=1$ if and only if $|\Gins(J_{\F(G)})|=1$.

Let $G$ be a graph on $[n]$.
If $G$ contains at most two edges,
then $|\Gins(J_G)|=1$ since there exists only one strongly stable ideal with the same Hilbert function as $J_G$.
An example of a graph $G$ with $|\Gins (J_G)|>1$ 
appears when $G$ has more than three edges.
Indeed,
it is not hard to see that the following graphs satisfy $|\Gins(J_G)|>1$
(Proposition \ref{threegraphs}).

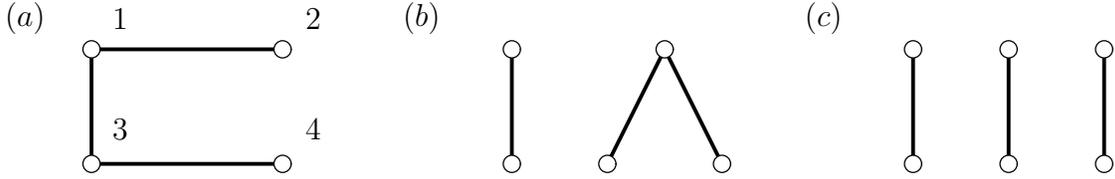
\vspace{10pt}

\begin{picture}(36.47, 22.63)(12.25, -30.63)
\put (11, -10){\makebox (0, 0){${(a)}$}}%
\put (47, -10){\makebox (0, 0){$1$}}%
\put (120, -10){\makebox (0, 0){$2$}}%
\put (47, -52){\makebox (0, 0){$3$}}%
\put (120, -52){\makebox (0, 0){$4$}}%

\put (161, -10){\makebox (0, 0){${(b)}$}}%
\put (313, -10){\makebox (0, 0){${(c)}$}}%

\special{pn 20}%
\special{pa 500 300}%
\special{pa 1500 300}%
\special{fp}%

\special{pn 20}%
\special{pa 500 900}%
\special{pa 1500 900}%
\special{fp}%

\special{pn 20}%
\special{pa 500 300}%
\special{pa 500 900}%
\special{fp}%

\special{pn 8}%
\special{sh 0.000}%
\special{ar 500 300 46 46  0.0000000 6.2831853}%

\special{pn 8}%
\special{sh 0.000}%
\special{ar 500 900 46 46  0.0000000 6.2831853}%

\special{pn 8}%
\special{sh 0.000}%
\special{ar 1500 300 46 46  0.0000000 6.2831853}%

\special{pn 8}%
\special{sh 0.000}%
\special{ar 1500 900 46 46  0.0000000 6.2831853}%


\special{pn 20}%
\special{pa 3500 300}%
\special{pa 3200 900}%
\special{fp}%

\special{pn 20}%
\special{pa 2700 300}%
\special{pa 2700 900}%
\special{fp}%

\special{pn 20}%
\special{pa 3500 300}%
\special{pa 3800 900}%
\special{fp}%

\special{pn 8}%
\special{sh 0.000}%
\special{ar 2700 900 46 46  0.0000000 6.2831853}%

\special{pn 8}%
\special{sh 0.000}%
\special{ar 2700 300 46 46  0.0000000 6.2831853}%

\special{pn 8}%
\special{sh 0.000}%
\special{ar 3200 900 46 46  0.0000000 6.2831853}%

\special{pn 8}%
\special{sh 0.000}%
\special{ar 3800 900 46 46  0.0000000 6.2831853}%

\special{pn 8}%
\special{sh 0.000}%
\special{ar 3500 300 46 46  0.0000000 6.2831853}%

\special{pn 20}%
\special{pa 4800 300}%
\special{pa 4800 900}%
\special{fp}%

\special{pn 20}%
\special{pa 5300 300}%
\special{pa 5300 900}%
\special{fp}%

\special{pn 20}%
\special{pa 5800 300}%
\special{pa 5800 900}%
\special{fp}%

\special{pn 8}%
\special{sh 0.000}%
\special{ar 4800 300 46 46  0.0000000 6.2831853}%

\special{pn 8}%
\special{sh 0.000}%
\special{ar 4800 900 46 46  0.0000000 6.2831853}%

\special{pn 8}%
\special{sh 0.000}%
\special{ar 5300 300 46 46  0.0000000 6.2831853}%

\special{pn 8}%
\special{sh 0.000}%
\special{ar 5300 900 46 46  0.0000000 6.2831853}%

\special{pn 8}%
\special{sh 0.000}%
\special{ar 5800 300 46 46  0.0000000 6.2831853}%

\special{pn 8}%
\special{sh 0.000}%
\special{ar 5800 900 46 46  0.0000000 6.2831853}%

\end{picture}%
\vspace{50pt}

Our result will explain that these three examples are essential if we characterize graphs $G$ with $|\Gins(J_G)|=1$
from the viewpoint of induced subgraphs.(See Theorem \ref{main1} (v) below.)

To state our result,
the notation about near cones is required.
Let $\Gamma$ be a simplicial complex on $[n]$.
We say that $\Gamma$ is a \textit{near cone} with respect to a vertex $v \in [n]$ if $\Gamma$ satisfies
$(S \setminus \{j\} ) \cup \{v\}$ for all $S \in \Gamma$ and $j \in S$.
If $\Gamma$ is a graph,
then the structure of a near cone is quite simple.
Let $G$ be a graph on $[n]$,
$G-v$ the graph obtained by deleting a vertex $v \in [n]$ from $G$
and $\deg_G(v)=|\{t: \{t,v\} \in G\}|$ the \textit{degree} of the vertex $v$ on $G$.
Then, $G$ is a near cone with respect to $v$ if and only if
$$G \supset  \{\{v\}\} \cup \{ \{v,t\}: t \in [n] \setminus \{v\} \mbox{ with }\deg_G(t)>0 \} \cup (G-v).$$
Thus if $G$ is a near cone with respect to a vertex $v\in [n]$,
then the structure of $G$ is determined from $G-v$ and $\deg_G(v)$.
We say that a graph $G$ is a $k$-\textit{near cone of $G'$}
if there exists a sequence $v_1,v_2,\dots,v_k \in [n]$ such that
$G- \{v_1,\dots,v_{t-1}\}$ is a near cone with respect to $v_t$ for $t=1,2,\dots,k$ and $G'=G-\{v_1,\dots,v_k\}$
is the graph on $[n] \setminus \{v_1,\dots,v_k\}$ obtained from $G$
by deleting vertices $v_1,\dots,v_k$.

We say that a graph $G$ is \textit{semi-complete} (resp.\ \textit{semi-complete bipartite})
if we can obtain a complete graph (resp.\ complete bipartite graph)
by deleting isolated vertices from $G$,
in other words,
$G$ is a union of a complete graph (resp.\ complete bipartite graph) and isolated vertices.
The main result of this paper is the following.

\begin{theorem} \label{main1}
Let $G$ be a graph on $[n]$.
The following conditions are equivalent.
\begin{itemize}
\item[(i)] $|\T (J_{\F(G)})|=1$;
\item[(ii)]  $|\Gins(J_{\F(G)})|=1$;
\item[(iii)] $|\Gins(J_G)|=1$;
\item[(iv)] $|\T(J_G)|=1$;
\item[(v)] $G$ and its complementary graph $\overline G$ contains none of the graphs $(a)$, $(b)$ and $(c)$ as an induced subgraph;
\item[(vi)] $G$ is a $k$-near cone of a semi-complete bipartite graph or of a disjoint union of two semi-complete graphs for some $k \geq 0$.
\end{itemize}
\end{theorem}

We also consider edge ideals.
Let $R=K[x_1,\dots,x_n]$ be the polynomial ring in $n$ variables over a field $K$
with char$(K)=0$.
For a homogeneous ideal $I \subset R$,
we write $\gin_\sigma (I)$ for the generic initial ideal with respect to a term order $\sigma$
and define 
$$\gins(I)=\{\gin_\sigma(I)\subset R:\sigma \mbox{ is a term order with }x_1>_\sigma \cdots>_\sigma x_n\}.$$

Let $G$ be a graph on $[n]$.
The \textit{edge ideal $I(G)\subset R$ of $G$} is the ideal generated by all squarefree monomials $x_ix_j$ with $\{i,j \} \in G$.
The second result of this paper is the following.

\begin{theorem}\label{main2}
Let $G$ be a graph on $[n]$. Assume that $\mathrm{char}(K)=0$.
Then we have $|\gins(I(G))|=1$ if and only if $G$ is a semi-complete bipartite graph.
\end{theorem}

Since the edge ideal $I(G)$ of $G$ is equal to the Stanley--Reisner ideal of the flag complex
of the complementary graph of $G$,
the above Theorem \ref{main2} determines all flag complexes
$\Gamma$ whose Stanley--Reisner ideal $I_{\Gamma}$ satisfies
$|\gins(I_{\Gamma})|=1$.

This paper is organized as follows:
In \S 2,
we study some properties for the homogeneous component of degree $2$ of generic initial ideals.
In \S 3,
we introduce some techniques which will be used to prove main theorems.
In \S 4,
we will show that the exterior face ideal $J_{\F(G)}$
of the flag complex $\F(G)$ of a semi-complete bipartite graph $G$
satisfies $|\Gins(J_{\F(G)})|=1$.
In \S 5,
the proof of Theorem \ref{main1} will be given.
In \S 6,
we will consider edge ideals $I(G)$ satisfying $|\gins(I(G))|=1$.

\section{The homogeneous component of degree $2$ of generic initial ideals}

In this section,
we study some properties for the homogeneous component of degree $2$ of strongly stable ideals.
Although we only consider the exterior algebra,
all of the lemmas in this section can be proved for homogeneous ideals in the polynomial ring
over a field $K$ with char$(K)=0$ in the same way.

First, we recall the fundamental theorem for generic initial ideals.
Let $K$ be an infinite field,
$V$ an $n$-dimensional $K$-vector space with basis $e_1,\dots,e_n$,
$E=\bigoplus_{k=0}^n \bigwedge^k V$ the exterior algebra of $V$
and $\GL$ the general linear group with coefficients in $K$.
Any $\varphi=(a_{ij}) \in \GL$ induces an automorphism of the graded $K$-algebra $E$
defined by
$$\varphi(e_i)= \sum_{k=1}^n a_{ki} e_k \ \ \ \ \mbox{ for } i=1,2,\dots,n.$$

\begin{lemma}[Galligo, Bayer--Stillman and Aramova--Herzog--Hibi] \label{gbs}
Let $J \subset E$ be a homogeneous ideal and $\sigma$ a term order.
Then, there exists a nonempty Zariski open subset $U \subset \GL$
such that $\init_{\sigma}( \varphi(J))$ is constant for all $\varphi \in U$.
Furthermore,
$U$ meets nontrivially the set of all upper triangular invertible matrices.
\end{lemma}

The above initial ideal $\init_\sigma ( \varphi (J))$ with $\varphi \in U$
is called the \textit{generic initial ideal of $J$ with respect to the term order $\sigma$},
and will be denoted $\Gin_\sigma (J)$.
We list some basic properties below
(see, e.g., \cite[Lemma 3.3]{BNT} for (i) -- (iv) and \cite[Corollary 2.3]{HM1} for (v)).

\begin{lemma} \label{itumono}
Let $J= \bigoplus_{d=0}^n J_d\subset E$ be a homogeneous ideal and $\sigma$ a term order
induced by $e_1> \cdots >e_n$,
where $J_d$ is the homogeneous component of degree $d$ of $J$.
Then
\begin{itemize} 
\item[(i)] $\Gin_\sigma (J)$ is strongly stable;
\item[(ii)]
if $J$ is strongly stable, then one has $\Gin_\sigma(J)=J$;
\item[(iii)] $J$ and $\Gin_\sigma(J)$ have the same Hilbert function,
that is, $\dim_K \Gin_\sigma(J)_d=\dim_K J_d$ for all integers $d \geq 0$;
\item[(iv)]
if $J \subset J'$ are homogeneous ideals in $E$,
then $\Gin_\sigma (J) \subset \Gin_\sigma(J')$;
\item[(v)] $\Gin_\sigma(J) = \Gin_\sigma( \varphi (J))$  for any $\varphi \in \GL$.
\end{itemize}
\end{lemma}

We also note similar properties for transformed strongly stable ideals.

\begin{lemma}\label{itumono2}
Let $J \subset E$ be a homogeneous ideal and $J' \in \T(J)$.
Then
\begin{itemize}
\item[(i)] $J$ and $J'$ have the same Hilbert function;
\item[(ii)]
if $J \subset I$ are homogeneous ideals in $E$,
then there exists an $I' \in \T(I)$ such that $J'\subset I'$;
\item[(iii)] $\T(J) = \T( \varphi (J))$  for any $\varphi \in \GL$.
\end{itemize}
\end{lemma}

\begin{proof}
The statement (i) follows from the fact that Hilbert functions do not change
by taking initial ideals.
The statement (ii) easily follows by using the fact that
if $J \subset I$ then $\init_\sigma (\varphi(J)) \subset \init_\sigma (\varphi(I))$ for any $\varphi \in \GL$
and any term order $\sigma$.
On the other hand, (iii) is obvious from the definition of transformed strongly stable ideals.
\end{proof}

Let $R=K[x_1,\dots,x_n]$ be the polynomial ring in $n$ variables with  each $\deg(x_i)=1$.
The \textit{graded Betti numbers} $\beta_{ij}(I)$ of a homogeneous ideal $I \subset R$ are the integers defined by
$$\beta_{ij}(I)= \dim_K \mathrm{Tor}_i (I,K)_j.$$
Studying graded Betti numbers is one of the current trend in computational commutative algebra.
Also, many nice relations between generic initial ideals and graded Betti numbers
have been discovered (see, e.g., \cite{AHH}, \cite{C}, \cite{C2}, \cite{H} and \cite{HM1}). 

Let $J$ be a monomial ideal in $E$.
We write $J^* \subset R$ for the ideal in the polynomial ring $R$
generated by all squarefree monomials $x_{i_1}x_{i_2} \cdots x_{i_k}$ with $e_{i_1}\wedge e_{i_2} \wedge \cdots \wedge e_{i_k} \in J$.
For an integer $1 \leq k \leq n$,
let
$$ \mathrm{min}_{ \leq k}(J,d)= |\{ e_S \in J_d: \min (S) \leq k\}|$$
and 
$$ \mathrm{max}_{ \leq k} (J,d)= |\{ e_S \in J_d: \max (S) \leq k\}|.$$

The homogeneous component of degree $2$ of strongly stable ideals has a simple structure.
Indeed, the next lemma easily follows from \cite[Corollary 3.6]{H}.

\begin{lemma}[{\cite[Lemma 3.1]{Mf}}] \label{equivalent}
Let $J \subset E$ and $J' \subset E$ be strongly stable ideals 
which do not contain any monomial of degree $1$.
The following conditions are equivalent.
\begin{itemize}
\item[(i)] $J_2 =J'_2$; 
\item[(ii)] $\max_{\leq k}( J,2) = \max_{\leq k}(J',2)$\ for all $k$;  
\item[(iii)] $\min_{\leq k}( J,2) = \min_{\leq k}(J',2)$\hspace{5pt} for all $k$;  
\item[(iv)] $\beta_{ii+2}( J^*) = \beta_{ii+2}( (J')^*)$\hspace{22pt} for all $i$.
\end{itemize}
\end{lemma}

Next, we recall the nice relation between generic initial ideals and transformed strongly stable ideals
found by Conca.
Let $\sigma$ be a term order and $e_S$ a monomial in $E$.
Set
$$ m_{\sigma,e_S} (J) = |\{e_T \in J: e_T \geq _\sigma e_S, \deg(e_T)= \deg(e_S)\}|.$$

\begin{lemma} \label{conca}
Let $\sigma$ be a term order and $J \subset E$ a homogeneous ideal.
If $J' \in \T (J)$, then, for any monomial $e_S \in E$,
one has
$$ m_{\sigma,e_S} (\Gin_\sigma(J)) \geq m_{\sigma,e_S}(J').$$
\end{lemma}

\begin{proof}
It follows from \cite[Proposition 2.4]{Mj} 
that, for any $e_S \in E$ and term order $\tau$,
one has
$m_{\sigma,e_S} (\Gin_\sigma(J)) \geq m_{\sigma,e_S}(\Gin_\sigma(\init_\tau(J)))$.
(Note that the same property for the case of the polynomial ring is \cite[Corollary 1.6]{C}.)
If $J'$ is a transformed strongly stable ideal of $J$ with
$J'= \init_{\sigma_k} ( \varphi_k ( \cdots \init_{\sigma_{2}} ( \varphi_2 (\init_{\sigma_1}(\varphi_1 (J)) )) \cdots )),$
then the above inequality together with Lemma \ref{itumono} (v) says that
\begin{eqnarray*}
 m_{\sigma,e_S} (\Gin_\sigma(J)) =m_{\sigma,e_S} (\Gin_\sigma(\varphi_1(J))) 
 &\geq&  m_{\sigma,e_S} (\Gin_\sigma(\init_{\sigma_1}( \varphi_1(J))) \\
 &\geq&  m_{\sigma,e_S} (\Gin_\sigma(\init_{\sigma_2}(\varphi_2(\init_{\sigma_1}( \varphi_1(J))))) \\
 &&\vdots\\
 &\geq&  m_{\sigma,e_S} (\Gin_\sigma(J')).
\end{eqnarray*}
Since $J'$ is strongly stable,
we have $\Gin_\sigma(J')=J'$.
Thus the assertion follows.
\end{proof}

Let $\lex$ (resp.\ $\rlex$) be the degree lexicographic
(resp.\ reverse lexicographic) order induced by $e_1> \cdots >e_n$.
We write $\Gin(J)$ for the generic initial ideal of $J$
with respect to $\rlex$.
The next lemma immediately follows from Lemma \ref{conca} together with the definition of $\lex$ and that of $\rlex$.

\begin{lemma}\label{minmax}
Let $J$ be a homogeneous ideal in $E$ and $J' \in \T(J)$.
Then
\begin{itemize}
\item[(i)] $\min_{\leq k}( \Gin_{\mathrm{lex}}(J),d) \geq \min_{\leq k}(J',d)$\ \  for all $k$ and $d$;
\item[(ii)] $\max_{\leq k}( \Gin(J),d) \geq \max_{\leq k}(J',d)$\ \ \hspace{3pt} for all $k$ and $d$.
\end{itemize}
\end{lemma} 

Moreover, the following fact is true.

\begin{proposition} \label{deg2}
Let $ J\subset E$ be a homogeneous ideal and $J' \in \T(J)$.
Then, for any integer $k \geq 1$, one has
\begin{eqnarray}
\mathrm{max}_{\leq k}( \Gin_{\mathrm{lex}}(J),2) \leq \mathrm{max}_{\leq k} (J',2) \leq \mathrm{max}_{ \leq k} ( \Gin(J),2). \label{parasol}
\end{eqnarray}
\end{proposition}

\begin{proof}
The right-hand side of (\ref{parasol})
is Lemma \ref{minmax} (ii).
We will show the left-hand side.
Fix an integer $k \geq 1$.
Let 
$$q= \max\{t:e_t \wedge e_k \in J'\}$$
where we let $q=0$ if $e_1 \wedge e_k  \not \in J'$.

If $q \geq k-1$,
then $J'$ contains all monomials $e_i \wedge e_j$ with $\max\{i,j\} \leq k$,
and therefore we have the inequality (\ref{parasol}).
Hence we assume that $q < k-1$.

Since $J'$ is strongly stable,
$e_{q+1} \wedge e_k \not \in J'$ implies that
any monomial $e_i \wedge e_j$ with $i > q$  and $j \geq k$
does not belong to $J'$.
Thus we have
\begin{eqnarray}
\{e_i \wedge e_j\in J': \max\{i,j\} >k\}
= \{ e_i \wedge e_j \in J': \max\{i,j\} >k,\ \min\{i,j\} \leq q\}. \label{sword}
\end{eqnarray}
Moreover, since $e_i \wedge e_j \in J'$ for all $i,j$ with $i \leq k$ and $j \leq q$, a routine computation says that
\begin{eqnarray}
&&|\{ e_i \wedge e_j \in J': \max\{i,j\} >k,\ \min\{i,j\} \leq q\}| \nonumber\\
&&\hspace{50pt}= \mathrm{min}_{ \leq q}(J',2) -\{(k-1)+(k-2)+ \cdots +(k-q)\}. \label{ball}
\end{eqnarray}
Also, in the same way as (\ref{ball}), a simple counting says that
\begin{eqnarray}
&&|\{e_i \wedge e_j\in \Gin_{\mathrm{lex}}(J): \max\{i,j\} >k\}| \nonumber\\ 
&&\hspace{50pt}\geq  |\{ e_i \wedge e_j \in \Gin_{\mathrm{lex}}(J): \max\{i,j\} >k,\ \min\{i,j\} \leq q\}| \nonumber\\
&&\hspace{50pt}\geq \mathrm{min}_{ \leq q}(\Gin_{\mathrm{lex}}(J),2) -\{(k-1)+(k-2)+ \cdots +(k-q)\}. \label{sleep}
\end{eqnarray}  
On the other hand,
for any strongly stable ideal $I \subset E$, one has
\begin{eqnarray}
|\{ e_i \wedge e_j \in  I: \max\{i,j\} >k\}|= \dim_K I_2 -\mathrm{max}_{\leq k}(I,2). \label{fire}
\end{eqnarray}
Then (\ref{sword}), (\ref{ball}), (\ref{sleep}) and (\ref{fire}) together with Lemma \ref{minmax} (i) say that
\begin{eqnarray*}
\dim_K J'_2 -\mathrm{max}_{\leq k} (J',2)
& = & \mathrm{min}_{ \leq q}(J',2) -\{(k-1)+ \cdots +(k-q)\}\\
& \leq & \mathrm{min}_{ \leq q}(\Gin_{\mathrm{lex}}(J),2) -\{(k-1)+ \cdots +(k-q)\}\\
& \leq & |\{e_i \wedge e_j\in \Gin_{\mathrm{lex}}(J): \max\{i,j\} >k\}|\\
&=& \dim_K \Gin_{\mathrm{lex}}(J)_2 -\mathrm{max}_{\leq k} (\Gin_{\mathrm{lex}}(J),2).
\end{eqnarray*}  
Since $J'$ and $\Gin_{\mathrm{lex}}(J)$ have the same Hilbert function, we have $\max_{\leq k}(J',2) \geq \max_{\leq k}( \Gin_{\mathrm{lex}}(J),2)$ as desired.
\end{proof}

It is known that the inequality (\ref{parasol})
yields the inequality of graded Betti numbers.
Indeed, the next corollary can be proved in the same way as
\cite[Proposition 3.6]{C2} and \cite[Theorem 4.4]{SQlex}.

\begin{cor} \label{betti}
Let $J \subset E$ be a homogeneous ideal and $J' \in \T (J)$.
Then one has
$$\beta_{ii+2}( \Gin_{\mathrm{lex}}(J)^*) \geq \beta_{ii+2}((J')^*) \geq \beta_{ii+2}( \Gin(J)^*) \ \ \mbox{ for all }i.$$
\end{cor}


Also, Corollary \ref{betti} together with Lemma \ref{equivalent} implies the following fact.
For an integer $d \geq 0$, define
$$ \Gins (J,d) = \{ J'_d: J' \in \Gins(J)\}$$
and
$$\T (J,d)= \{ J_d': J' \in \T(J) \}.$$

\begin{cor} \label{unique}
Let $J\subset E$ be a homogeneous ideal.
The following conditions are equivalent.
\begin{itemize}
\item[(i)] $\Gin_{\mathrm{lex}}(J)_2= \Gin(J)_2$;
\item[(ii)] $|\T(J,2)|=1$;
\item[(iii)] $|\Gins(J,2)|=1$.
\end{itemize}
\end{cor}

If $G$ is a graph,
then $J_G$ contains all monomials of degree $\geq 3$.
Thus $|\Gins(J_G,d)|=1$ for $d \geq 3$.
Also, $\Gins(J_G,d)=\{0\}$ if $d=1$ or $d=2$.
Then, by Corollary \ref{unique},
we have $|\Gins(J_G)=1|$ if and only if $\Gin_{\mathrm{lex}}(J_G)=\Gin(J_G)$.

\begin{exam}
Proposition \ref{deg2} and Corollaries \ref{betti} and \ref{unique} are false for a degree $d \geq 3$.
In fact, let $J$ be the strongly stable ideal generated by
\begin{eqnarray*}
\left\{
\begin{array}{r}
e_1\wedge e_2 \wedge e_3,e_1\wedge e_2 \wedge e_4,e_1\wedge e_2 \wedge e_5,e_1\wedge e_2 \wedge e_6,\\
e_1\wedge e_3 \wedge e_4,e_1\wedge e_3 \wedge e_5,e_1\wedge e_3 \wedge e_6,\\
e_2\wedge e_3 \wedge e_4,e_2\wedge e_3 \wedge e_5\hspace{65pt}
\end{array}
\right\}
\end{eqnarray*}
together with all monomials in $E$ of degree $4$.
Let 
$$J' =J+ (e_4 \wedge e_5 \wedge e_6).$$

Since $J$ is strongly stable, we have $\Gin_\sigma(J)=J$ for any term order $\sigma$.
Thus, we have $\Gin_\sigma(J') \supset \Gin_\sigma(J)=J$ by Lemma \ref{itumono} (iv).
On the other hand,
if $I \subset E$ is a strongly stable ideal which contains $J$ and has the same Hilbert function as $J'$,
then $I$ is either $J+(e_1 \wedge e_4 \wedge e_5)$ or $J+(e_2 \wedge e_3 \wedge e_6)$.
In fact, we have $\Gin(J')=\Gin_{\mathrm{lex}}(J')=J+(e_1 \wedge e_4 \wedge e_5)$.
However, if $\sigma$ is a term order induced by the weight $(10,9,8,3,2,1)$,
then we have $\Gin_\sigma(J')=J+(e_2 \wedge e_3 \wedge e_6)$.
Moreover, we have
$$\beta_{4\ 4+3}({\Gin_{\mathrm{lex}}(J')}^*)=2 <3=\beta_{4\ 4+3}({\Gin_{\sigma}(J')}^*).$$
Thus this example shows that Proposition \ref{deg2} and Corollaries \ref{betti} and \ref{unique}
are false for homogeneous components of degree $ \geq 3$.
\end{exam}

\begin{exam} \label{rei}
We will explain the advantage of considering $\T (J)$.
In general, it is not easy to determine $\Gin_\sigma(J)$
since we must compute $\init_\sigma (\varphi (J))$ for a generic matrix $\varphi \in \GL$.
Thus determining $\Gins(J)$ is quite difficult.
However, we do not need to consider a generic matrix $\varphi \in \GL$ when we consider transformed strongly stable ideals,
and $|\T(J,2)|>1$ implies $|\Gins(J)|>1$ by Corollary \ref{unique}.

Moreover, the following idea yields many transformed
strongly stable ideals by considering only elementary matrices in $ \GL$.

Let $J \subset E$ be a homogeneous ideal and $\sigma$ a term order induced by $e_1> \cdots >e_n$.
For positive integers $1 \leq a<b \leq n$,
write $\varphi_{a,b} \in \GL$ for the elementary matrix defined by
$\varphi_{a,b}(e_k)=e_k$ if $k \ne b$,
and $\varphi_{a,b}(e_b)=e_a +e_b$.
Then, it is not hard to show that
there exists a sequence of pairs of positive integers $(a_1,a_2),\dots,(a_p,b_p)$, where $a_k<b_k$ for each $k$,
such that
$$J'= \init_\sigma(\varphi_{a_p,b_p} ( \cdots \init_\sigma ( \varphi_{a_1,b_1}(J)) \cdots )$$
is strongly stable.
This monomial ideal $J'$ is a transformed strongly stable ideal of $J$.
In particular, if $J$ is a monomial ideal,
the above construction does not depend on term orders and is known as combinatorial shifting.
(See \cite[Lemma 8.3]{H}.)

For example, let $J=(e_1\wedge e_2, e_1\wedge e_3,e_3\wedge e_4)$.
Then, 
$$J'= \init_\sigma (\varphi_{1,3}(J))=(e_1\wedge e_2, e_1\wedge e_3,e_1\wedge e_4,e_2\wedge e_3 \wedge e_4)$$
and
$$J''=\init_\sigma (\varphi_{2,4}(J))=(e_1\wedge e_2, e_1\wedge e_3,e_2\wedge e_3)$$
are transformed strongly stable ideals of $J$.
Thus we have $|\T(J,2)|>1$, and therefore we have $|\Gins(J)|>1$ by Corollary \ref{unique}.
\end{exam}

\section{Some basic techniques}

In this section,
we will introduce two techniques which will be used for the proof of Theorems \ref{main1} and \ref{main2}.
The first one (Lemmas \ref{symcomplement} and \ref{complement}) is required
to prove $|\Gins(J_G)|=1$ and $|\gins(I(G))|=1$,
and the second one (Lemmas \ref{syminduced} and \ref{induced}) is required to prove $|\Gins(J_G)|>1$ and $|\gins(I(G))|>1$.

Let $K$ be an infinite field
and $R=K[x_1,\dots,x_n]$ the polynomial ring in $n$ variables with each $\deg (x_i)=1$.
A monomial ideal $I \subset R$ is called \textit{strongly stable} if
$u x_q \in I$ implies $u x_p \in I$ for all $1 \leq p < q \leq n$.
If $\mathrm{char}(K)=0$, then the generic initial ideal
$\gin_\sigma (I)$ of $I \subset R$ is strongly stable
for an arbitrary term order $\sigma$.
Thus we assume $\mathrm{char}(K)=0$ when we consider generic initial ideals in the polynomial ring $R$.
For a $K$-vector subspace $W \subset R_d$ and $\varphi \in \GL$,
let $\varphi(W)=\{\varphi(f):f \in W\}$ and $\init_\sigma(W)$ the $K$-vector space spanned by $\{ \init_\sigma(f):f \in W\}$.
Also, we define $\gin_\sigma (W)$, $\gins(W)$ and $\T (W)$ in the same way as \S 1.

Let $W \subset R_d$ be a $K$-vector space spanned by monomials.
We write $\overline W$ for the $K$-vector space spanned by all monomials $u\in R_d$ with $u \not \in V$.
For a term order $\sigma$,
we write $\sigma^{-1}$ for the term order defined by
$u >_{\sigma^{-1}} v$ if $\deg(u)> \deg (v)$ or $\deg(u) = \deg (v)$ and $ u<_\sigma v$.
The following fact was written in \cite[p.\ 133]{K} without a proof.

\begin{lemma} \label{symcomplement}
Let $W \subset R_2$ be a $K$-vector space spanned by monomials of
degree $2$ and $\sigma$ a term order.
If $\mathrm{char(K)}=0$
then
$$ \overline {\gin_\sigma (W)} = \gin_{\sigma^{-1}} ( \overline W).$$
\end{lemma}

\begin{proof}
Recall that if $U \subset \GL$ and $U' \subset \GL$ are nonempty Zariski open subsets
then $U \cap U'$ and $U^{-1}=\{ \varphi^{-1}: \varphi \in U\}$
are also nonempty Zariski open subsets of $\GL$.
This fact together with
Lemma \ref{gbs} says that
there exists a $\varphi=(a_{ij}) \in \GL$ such that
$\init_\sigma (\varphi(W))= \gin_\sigma(W)$ and $\init _{\sigma^{-1}} ( \varphi^{-1}( \overline W))=
 \gin_{\sigma^{-1}}( \overline W)$.
For each monomial $x_ix_j\in V$,
$\varphi(x_ix_j)$ can be written in the form
$$
\varphi(x_ix_j)= \sum_{1 \leq p< q \leq n} ( a_{pi}a_{qj}+a_{pj}a_{qi}) x_p x_q + \sum_{p=1,2,\dots,n} a_{pi}a_{pj}x_p^2.
$$
Let
\begin{eqnarray*}
\alpha_{(x_ix_j,x_px_q)} =
\left\{
\begin{array}{ll}
a_{pi}a_{qj}+a_{pj}a_{qi}, \ \ & \mbox{ if } p \ne q,\\
a_{pi}a_{pj}, \hspace{55pt} & \mbox{ if } p=q.
\end{array}
\right.
\end{eqnarray*}
Consider the  $\dim_K W \times { n+1 \choose 2}$ matrix 
$$X=( \alpha_{(x_ix_j,x_px_q)} )_{x_ix_j \in W,\ x_px_q \in R_2}.$$
For each monomial $x_px_q\in R_2$,
let $X_{(\sigma,x_px_q)}$ be the submatrix of $X$ which consists of the columns of $X$ indexed by
$x_sx_t$ with $x_sx_t \geq _\sigma x_p x_q$
and 
$X'_{(\sigma,x_px_q)}$ the submatrix of $X_{(\sigma,x_px_q)}$
which is obtained by removing the columns of $X_{(\sigma,x_px_q)}$ indexed by $x_px_q$.
Then, it is not hard to see that (see, e.g., \cite[Lemma 2.1]{HM1})
\begin{eqnarray}
x_px_q \not \in \init_\sigma(\varphi(W)) \ \ \ 
\mbox{ if and only if } \ \ \ \mathrm{rank}(X_{(\sigma,x_px_q)})= \mathrm{rank}( X'_{(\sigma,x_px_q)}). \label{LRM1}
\end{eqnarray}\smallskip

Next, consider the quotient space $A=R_2 / \varphi^{-1}(\overline W)$.
For any polynomial $f \in R_2$,
we write $[f]$ for its image in $A$.
Then, for any monomial $x_px_q\in R_2$,
the next fact follows from the definition of initial monomials.\begin{eqnarray}
x_px_q \in \init_{\sigma^{-1}}(\varphi^{-1}(\overline W)) \ \ \mbox{ if and only if } \ \ 
[x_px_q] \in \mathrm{span}\{
[x_sx_t] \in A: x_sx_t <_{\sigma^{-1}}x_px_q \},
\label{LRM2}
\end{eqnarray}
where $\mathrm{span}(V)$ with $V \subset A$
is the $K$-vector space spanned by all elements in $V$.
Also, since  $[\varphi^{-1}(x_sx_t)]=0$ for all $x_sx_t \not \in W$,
for any monomial $x_px_q$, we have
\begin{eqnarray}
[x_px_q] &=& [\varphi^{-1} ( \varphi(x_px_q))] \nonumber \\
&=& \bigg [ \varphi^{-1} \Big\{ \sum_{1 \leq i < j \leq n}(a_{pi}a_{qj}+a_{pj}a_{qi})x_ix_j+ \sum_{i=1,2,\dots,n} (a_{pi}a_{qi}) x_i^2\Big\} \bigg ] \nonumber \\
&=& \sum_{ {1 \leq i< j \leq n} \atop {x_i x_j \in  W}} (a_{pi}a_{qj}+a_{pj}a_{qi})\ [\varphi^{-1}(x_ix_j)]
+ \sum_{ {i=1,2,\dots,n} \atop {x_i^2 \in  W}} (a_{pi}a_{qi})\ [\varphi^{-1} (x_i^2)]. \label{tori}
\end{eqnarray}
Let 
\begin{eqnarray*}
\beta_{(x_ix_j,x_px_q)} =
\left\{
\begin{array}{ll}
a_{pi}a_{qj}+a_{pj}a_{qi}, \ \ \ &\mbox{ if } i \ne j,\\
a_{pi}a_{qi}, \hspace{54pt} &\mbox{ if } i=j.
\end{array}
\right.
\end{eqnarray*}
Consider the $\dim_K W \times {n+1 \choose 2}$ matrix
$$Y =( \beta_{(x_ix_j,x_px_q)})_{x_ix_j \in W,x_px_q \in R_2}.$$
Let $Y_{(\sigma^{-1},x_px_q)}$
be the submatrix of $Y$ which consists of the columns of $Y$ indexed by
$x_sx_t$ with $x_sx_t \leq _{\sigma^{-1}} x_p x_q$
and define $Y'_{(\sigma^{-1},x_px_q)}$ in the same way as $X'_{(\sigma,x_px_q)}$.
Since $\{[\varphi^{-1}(x_sx_t)] \in A: x_s x_t \in W\}$ is a $K$-basis of $A$,
(\ref{LRM2}) and (\ref{tori}) say that
\begin{eqnarray}
x_px_q \in \init_{\sigma^{-1}} ( \varphi^{-1}(\overline W))\ \ \  \mbox{ if and only if } \ \ \
\mathrm{rank}( Y_{(\sigma^{-1},x_px_q)}) = \mathrm{rank}( Y'_{(\sigma^{-1},x_px_q)}). \label{LRM3}
\end{eqnarray}
We will show $\mathrm{rank}(X_{(\sigma,x_px_q)})= \mathrm{rank}(Y_{(\sigma^{-1},x_px_q)})$ for any monomial $x_px_q$.
\smallskip

If we divide each $x_px_q$-th column vector, with $p \ne q$,
of $X$ by $2$,
then each $(x_ix_j,x_px_q)$-th entry becomes $\frac 1 2 ( a_{pi}a_{qj} + a_{pj}a_{qi})$.
On the other hand,
if we divide each $x_{i}x_j$-th low vector, with $i\ne j$, of $Y$ by $2$,
then each $(x_ix_j,x_px_q)$-th entry becomes $\frac 1 2 ( a_{pi}a_{qj}+a_{pj}a_{qi})$.
Thus the definition of $\sigma^{-1}$ says that
$\mathrm{rank}(X_{(\sigma,x_px_q)})= \mathrm{rank}(Y_{(\sigma^{-1},x_px_q)})$
for any monomial $x_px_q$.
Also, this fact says that
$\mathrm{rank}(X'_{(\sigma,x_px_q)})= \mathrm{rank}(Y'_{(\sigma^{-1},x_px_q)})$
for any monomial $x_px_q$.

Recall that $\init_\sigma (\varphi(W))=\gin_\sigma(W)$ and
$\init_{\sigma^{-1}} ( \varphi^{-1}(\overline W)) = \gin_{\sigma^{-1}} ( \overline W)$.
Then
(\ref{LRM1}) and (\ref{LRM3}) say that,
for any monomial $x_px_q$, we have
$x_px_q \in \overline {\gin_\sigma(W)}$ if and only if $x_px_q \in \gin_{\sigma^{-1}}( \overline W)$,
as desired.
\end{proof}

\begin{rem}
Lemma \ref{symcomplement} is false if char$(K) \ne 0$.
Let $W \subset R_2$ be the $K$-vector space spanned by $x_1^2$ and $x_2^2$.
If char$(K)=2$, then we have $\gin_\sigma(W)=W$ but
$\gin_\sigma(\overline W) \ne \overline W$  for any term order $\sigma$.
\end{rem}

A similar property is true in the exterior algebra.
(But we do not need to assume that char$(K)=0$.)
Let $W \subset E_2$ be a $K$-vector space spanned by monomials of degree $2$
and $\sigma$ a term order.
Define the $K$-vector space $\Gin_\sigma(W) \subset E_2$ and $\overline W \subset E_2$
in the same way as in the polynomial ring.

\begin{lemma} \label{complement}
Let $W \subset E_2$ be a $K$-vector space spanned by monomials of degree $2$
and $\sigma$ a term order.
Then 
$$ \overline {\Gin_\sigma(W)} =  \Gin_{\sigma ^{-1}} (\overline W).$$
\end{lemma}
\begin{proof}
We sketch the proof since it is similar to
the proof of Lemma \ref{symcomplement}.
Let $\varphi=(a_{ij}) \in \GL$ be a matrix
with
$\init_\sigma(\varphi(W))=\Gin_\sigma(W)$ and $\init_{\sigma^{-1}}(\varphi( \overline W))=\Gin_{\sigma^{-1}}(\overline W)$.
Then, for each monomial $e_i \wedge e_j \in E$,
we have
\begin{eqnarray}
&&\varphi(e_i \wedge e_j)=
\sum_{1 \leq p<q \leq n}(a_{pi}a_{qj}-a_{pj}a_{qi})e_p \wedge e_q \label{P1}.
\end{eqnarray}

For any monomial $e_p \wedge e_q\in E_2$,
we write $[e_p \wedge e_q]$ for its image in $E_2/ \varphi^{-1}(\overline W)$.
Then,
for any monomial $e_p \wedge e_q\in E_2$,
we have
\begin{eqnarray}
&&[e_p \wedge e_q] = \sum_{ {1 \leq i< j  \leq n},\ {e_i \wedge e_j\in W}}
(a_{pi}a_{qj}-a_{pj}a_{qi})[\varphi^{-1} (e_i \wedge e_j) ]. \label{P2}
\end{eqnarray}
By using (\ref{P1}) and (\ref{P2}),
the statement follows in the same way as Lemma \ref{symcomplement}
\end{proof}

Lemmas \ref{symcomplement} and \ref{complement} would be true for an arbitrary degree.
However, we only require these facts for a degree $d=2$ in this paper.

\begin{exam}  \label{reisann}
By using Lemma \ref{complement},
we can compute $\Gin_{\mathrm{lex}}(\overline W)$ from $\Gin(W)$
if $W\subset E_2$ is a $K$-vector subspace spanned by monomials.
Let $W \subset E_2$ be the $K$-vector space spanned by $\{ e_1 \wedge e_2, e_2 \wedge e_3, e_3 \wedge e_4\}$.
Then $\Gin(W)$ is spanned by $\{ e_1 \wedge e_2, e_1 \wedge e_3, e_2 \wedge e_3\}$.
Thus, by Lemma \ref{complement},
$\Gin_{\mathrm{rev}^{-1}}(\overline W)$ is spanned by $\{ e_1 \wedge e_4, e_2 \wedge e_4, e_3 \wedge e_4\}$.
Recall that $<_{\mathrm{rev}^{-1}}$ is the lexicographic order induced by $e_1<e_2<\cdots <e_n$.
This fact says that $\Gin_{\mathrm{lex}}(\overline W)$ is spanned by $\{ e_1 \wedge e_2, e_1 \wedge e_3, e_1 \wedge e_4\}$.
\end{exam}

Next,
we will introduce the property which will be used to prove $|\Gins(J_G)|>1$ and $|\gins(I(G))|>1$.
We first require the following simple fact.

\begin{lemma} \label{partial}
Let $W\subset R_2$ be a $K$-vector space spanned by monomials.
If $x_sx_t \not \in W$ for all $(s,t)$ with $s\geq p$ and $t \geq q$,
then $x_px_q \not \in \gin_\sigma(W)$
for any term order $\sigma$ with $x_1>_\sigma x_2 >_\sigma \cdots  >_\sigma x_n$.
\end{lemma}

\begin{proof}
Lemma \ref{gbs} says that there exists an upper triangular matrix $\varphi \in \GL$ such that
$\init_\sigma (\varphi(W)) = \gin_\sigma(W)$.
Then, by the assumption, $x_px_q$ does not appear in $\varphi(u)$ for any monomial $u$ in  $W$.
Thus we have $x_px_q \not \in \init_\sigma (\varphi(W))=\gin_\sigma(W)$.
\end{proof}

Let $W \subset R_d$ be a $K$-vector space spanned by monomials.
For a subset $S=\{i_1,i_2,\dots,i_k\} \subset [n]$,
we write
$$ W_S=\{f \in W: f \in K[x_{i_1},x_{i_2},\dots,x_{i_k}]\}.$$

\begin{lemma} \label{syminduced}
Assume that $\mathrm{char}(K)=0$.
Let $W \subset R_2$ be a $K$-vector space spanned by monomials of degree $2$.
If $|\T (W)|=1$, then we have $|\T(W_S)|=1$
for any subset $S \subset [n]$.
\end{lemma}

\begin{proof}
It is enough to prove the claim when $|S|=n-1$.
We may assume that $S=\{1,2,\dots,n-1\}$ by Lemma \ref{itumono} (v).
For $\varphi \in \mathrm{GL}_{n-1}(K)$, define $\tilde \varphi \in \GL$ by $\tilde \varphi (x_i)=\varphi(x_i) $ for $i=1,2,\dots,n-1$
and $\tilde \varphi (x_n)=x_n$.
Let $\sigma $ be a term order and $t=|\{i:x_ix_n \in W\}|$.
Then,
in the same way as the proof of Lemma \ref{symcomplement},
 there exists $\varphi \in \mathrm{GL}_{n-1}(K)$ such that 
\begin{eqnarray}
\init_\sigma (\tilde \varphi (W)) =
\left\{
\begin{array}{ll}
 \gin_\sigma(W_S) +\mathrm{span} \{x_1x_n,x_2x_n,\dots,x_tx_n\}, \hspace{29pt}\ \ \mbox{if } x_n^2 \not \in W, \\
 \gin_\sigma(W_S) +\mathrm{span} \{x_1x_n,x_2x_n,\dots,x_{t-1}x_n,x_n^2\}, \ \ \mbox{if } x_n^2 \in W
 \end{array}
\right. \label{MG}
\end{eqnarray}
for $\sigma=\lex$ and $\sigma=\rlex$.
By Corollary \ref{unique},
if $|\T(W_S)|>1$ then $\gin(W_S) \ne \gin_{\mathrm{lex}} (W_S)$.
We will show that $\gin(W_S) \ne \gin_{\mathrm{lex}} (W_S)$ implies $|\T(W)|>1$.

Suppose that $\gin(W_S) \ne \gin_{\mathrm{lex}}(W_S)$.
If $t=0$, then $|\T(W)|>1$ is obvious.
Thus we assume $t >0$.
Let $x_px_q$ be the lexicographically largest monomial in $\gin(W_S) \setminus \gin_{\mathrm{lex}}(W_S)$,
where $1 \leq p \leq q \leq n$.
Then we have $x_{p-1}x_q \in \gin_{\mathrm{lex}}(W_S)$.
Set 
$$W_1 = \init_{\mathrm{rev}} ( \tilde \varphi (W))\ \  \mbox{ and }\ \  W_2 = \init_{\mathrm{lex}}( \tilde \varphi(W)).$$

\noindent
\textbf{[Case 1]}
Assume that $t<p$.
Then, by (\ref{MG}), $\init_{\mathrm{lex}} (\varphi_{t,n}(W_2))$ does not contain any monomial $x_l x_m$ with $l \geq p$ and $m \geq q$,
where $\varphi_{t,n} \in \GL$ is the matrix defined in Example \ref{rei}.
Thus Lemma \ref{partial} says that
$x_px_{q}\not \in \gin(\init_{\mathrm{lex}} (\varphi_{t,n}(W_2))) \in \T(W)$.
On the other hand,
$W_1$ contains all monomials $x_lx_m$ with $l \leq p$ and $m \leq q$.
Thus, by Lemma \ref{itumono} (ii) and (iii),
we have $x_px_{q} \in \gin( W_1) \in \T(W)$.
Hence we have $|\T (W)| >1$.
\medskip

\noindent
\textbf{[Case 2]}
Assume that $p\leq t \leq q$.
Since $x_{p-1}x_q \in \gin_{\mathrm{lex}}(W_S)$, we have
\begin{eqnarray*}
&&\init_{\mathrm{lex}} ( \varphi_{q,n}(W_2))\\
&&=
\left\{
\begin{array}{lll}
\gin_{\mathrm{lex}} (W_S) +
\mathrm{span}\{ x_1x_n,\dots,x_{p-1}x_n,x_px_q,\dots,x_tx_q\} \ \ \mbox{if } x_n^2 \not\in W,\\
\left.
\begin{array}{l}
\hspace{-6pt} \gin_{\mathrm{lex}} (W_S) \\
\ \ \ \ \ + \mathrm{span}\{ x_1x_n,\dots,x_{p-1}x_n,x_px_q,\dots, x_{t-1}x_q,x_q^2\}
\end{array}
\right.
\ \hspace{1pt} \mbox{if } x_n^2 \in W\mbox{ and } t>p,\\
\gin_{\mathrm{lex}} (W_S) +
\mathrm{span}\{ x_1x_n,\dots,x_{p-1}x_n,x_{q^2}\} \ \hspace{58.5pt} \mbox{if } x_n^2 \in W \mbox{ and } t=p.\\
\end{array}
\right.
\end{eqnarray*}
Then Lemma \ref{partial} says that
$x_px_{q+1}\not \in\gin( \init_{\mathrm{lex}}( \varphi_{q,n}(W_2)))$.
Set
$$W_3=\init_{\mathrm{lex}}(\varphi_{p,q+1}(\init_{\mathrm{lex}} ( \varphi_{q+1,n} (W_1)))).$$
Then, by (\ref{MG}),
$W_3$ contains all monomials $x_lx_m$ with $l \leq p$ and $m \leq q+1$.
Thus, by Lemma \ref{itumono} (ii) and (iii),
we have $x_px_{q+1} \in \gin( W_3)$.
Hence we have $|\T (W)| >1$.
\smallskip

\noindent
\textbf{[Case 3]}
Assume that $t>q$.
We substitute $x_0$ for $x_n$ and regard $W_1$ and $W_2$ as subsets of $K[x_0,x_1,\dots,x_{n-1}]$.
Let $\tau$ be a term order satisfying $x_0>_\tau x_1>_\tau \cdots >_\tau x_{n-1}$.
Then, by Lemma \ref{partial},
we have $x_px_q \not \in \gin_\tau(W_2)$.
On the other hand, (\ref{MG}) says that
$\init_{\tau}( \varphi_{0,t}(W_1))$
contains all monomials $x_lx_m$ with $0 \leq l \leq p$ and $0 \leq m \leq q$.
Thus we have $x_px_q \in \gin_\tau (\init_{\tau}( \varphi_{0,t}(W_1)))$ by Lemma \ref{itumono} (ii) and (iii).

We regard $\gin_\tau(W_2)$ and $\gin_\tau (\init_{\tau}( \varphi_{0,t}(W_1)))$ as subsets of $K[x_1,\dots,x_n]$
again by substituting $x_n$ for $x_0$.
Consider the permutation $\pi \in \GL$ with $\pi(x_n)=x_1$ and $\pi(x_i)=x_{i+1}$ for $i=1,2,\dots,n-1$.
Then,
$x_{p+1}x_{q+1} \not \in \pi(\gin_\tau(W_2)) \in \T(W)$ and $x_{p+1}x_{q+1} \in \pi(\gin_\tau (\init_{\tau}( \varphi_{0,t}(W_1))))
\in \T(W)$.
Hence we have $|\T(W)|>1$.
\end{proof}
Lemma \ref{syminduced} immediately implies the next corollary.

\begin{cor} \label{edgeinduced}
Let $G$ be a graph on $[n]$.
If $|\T(I(G),2)|=1$, then, for any subset $S \subset [n]$,
one has $|\T(I(G_S),2)|=1$,
where $G_S$ is the induced subgraph of $G$ on $S$.
\end{cor}

Also, the following fact can be proved in the same way as Lemma \ref{syminduced}.

\begin{lemma} \label{induced}
Let $G$ be a graph on $[n]$.
If $|\T(J_G)|=1$ then, for any subset $S \subset [n]$,
one has $|\T(J_{G_S})|=1$.
\end{lemma}

\begin{proof}
Let $W=(J_G)_2$ and $W_S=(J_{G_S})_2$.
It is enough to prove that
if $|\T(W_S)|>1$ then $|\T(W)|>1$.
Set $t =|\{i\in [n] \setminus \{t\}:e_i \wedge e_t: \in J_G\}|$.
Then
there exists a matrix $\varphi \in \GL$ such that
$$\init_\sigma (\varphi(W))=\gin_\sigma(W_S)+
\mathrm{span}\{e_n \wedge e_1,\dots,e_n \wedge e_t\}$$
for $\sigma=\lex$ and $\sigma=\rlex$.
What we must prove is that
if $\Gin_{\mathrm{lex}}(W_S) \ne \Gin(W_S)$
then $|\T(W)|>1$.
Set $W_1 = \init_{\mathrm{rev}}(\varphi(W))$ and $W_2= \init_{\mathrm{lex}}(\varphi(W))$.
Assume that $\Gin(W_S) \ne \Gin_{\mathrm{lex}}(W_S)$
and $e_p \wedge e_q$ is the lexicographically largest monomial
in $\Gin(W_S) \setminus \Gin_{\mathrm{lex}}(W_S)$, where $1 \leq p <q \leq n$.

If $t< q$,
then we have $e_p \wedge e_q \in \Gin(W_1)$ and
$e_p \wedge e_q \not \in \Gin(W_2)$ in the same way as the [Case 1]
in Lemma \ref{syminduced}.
On the other hand,
if $t \geq q$,
then we have $e_{p+1} \wedge e_{q+1} \in \Gin( \pi (W_1) )$ and
$e_{p+1} \wedge e_{q+1} \not \in \Gin( \pi (W_2) )$
in the same way as the [Case 3] in Lemma \ref{syminduced},
where $\pi$ is the permutation defined by $\pi(e_n)=e_1$ and
$\pi(e_k)=e_{k+1}$ for $k=1,\dots,n-1$.
In both cases, we have $|\T(W)|>1$ as desired.
\end{proof}

\section{complete bipartite graphs and generic initial ideals}

In this section,
we will show (vi) $\Rightarrow$ (i) of Theorem \ref{main1}.

Let $\Gamma$ be a simplicial complex on $[n]$.
Thus ${\Gamma}$
is a collection of subsets of $[n]$ such that
(i) $\{ j \} \in {\Gamma}$ for all $j \in [n]$ and
(ii) if  $T \in {\Gamma}$ and  $S \subset T$ then $S \in {\Gamma}$.
The \textit{dimension} of $\Gamma$ is $\dim \Gamma = \max\{|S|:S \in \Gamma\}-1$.
For an integer $k \geq 1$,
we write $\Gamma_{k-1}=\{S \in \Gamma: |S|=k\}$.
The \textit{cone $\{n+1\} * \Gamma $ of $\Gamma$ over $n+1$} is the simplicial complex on $[n+1]$ generated by $\{ \{n+1\}\cup S: S \in \Gamma\}$.

Let $\sigma$ be a term order.
For any simplicial complex $\Gamma$,
define the simplicial complex $\Delta^\sigma ( \Gamma)$ on the same vertex set as $\Gamma$ by
$$J_{\Delta^\sigma (\Gamma)}= \Gin_\sigma(J_\Gamma).$$ 
In particular,
we write $\dele \Gamma$ for the simplicial complex defined by $J_{\dele \Gamma} = \Gin (J_\Gamma)$,
and call $\dele \Gamma$ the 
\textit{exterior algebraic shifted complex} of $\Gamma$.
This construction says that
knowing $\Delta^\sigma  (\Gamma)$ is equivalent to knowing $\Gin_\sigma(J_\Gamma)$.

Relations between cones and exterior algebraic shifted complexes
are well studied in \cite[\S 5]{N}.
In particular,
the following nice relation is known.

\begin{lemma}[Kalai {\cite[\S 2]{K}}]\label{cone}
For any simplicial complex $\Gamma$ on $[n]$,
one has
$$\dele{\{n+1\} * \Gamma}= \{n+1\} * \dele{ (\Gamma)}.$$
\end{lemma}

Let $\Gamma$ and $\Gamma'$ be simplicial complexes on $[n]$.
We say that $\Gamma'$ is combinatorially isomorphic to $\Gamma$
if there exists a permutation $\pi:[n] \to [n]$ such that
$\Gamma'=\{\pi(S): S \in \Gamma\}$.
We write $\Gamma \cup \{n+1\}$ for the simplicial complex on $[n+1]$ obtained by
adding the vertex $n+1$ to $\Gamma$.

\begin{lemma} \label{transcone}
Let $\Gamma$ be a $(d-1)$-dimensional simplicial complex on $[n-1]$.
If $|\T (J_{\Gamma\cup\{n\}})|=1$,
then one has $|\T (J_{\{n\} * \Gamma})|=1$.
\end{lemma}

\begin{proof}
We use induction on $d$.
If $d=1$, then $\{n\}*\Gamma$ is isomorphic to a shifted graph,
that is,
there exists a permutation $\pi$ such that $\pi (J_{\{n\}*\Gamma})$ is strongly stable.
Thus, by Lemma \ref{itumono} (ii) and Corollary \ref{unique}, we have  $|\T (J_{\{n\} * \Gamma})|=|\Gins (J_{\{n\} * \Gamma})|=1$.

Assume that $d>1$.
Let $J' \in \T (J_{\{n\} * \Gamma})$
and $\Gamma'$ the simplicial complex defined by $J_{\Gamma'}=J'$.
Set $\Sigma=\bigcup_{k=0}^{d-2}\Gamma_k$.
Then $\Sigma$ is a $(d-2)$-dimensional simplicial complex on $[n-1]$.
Also, we have $|\T (J_{ \Sigma \cup \{n\}},k)|=|\T(J_{\Gamma \cup \{n\}},k)|=1$ for $k\leq d-1$ by the assumption
and $|\T (J_{ \Sigma \cup \{n\}},k)|=1$ for $k \geq d$ since $J_{ \Sigma \cup \{n\}}$ contains all monomials in $E$ of degree $\geq d$.
Thus $|\T (J_{\{n\} * \Sigma})|=1$ by the induction hypothesis.
Since $\{n\}* \Gamma \supset \{n\}* \Sigma$,
we have $J_{\{n\}* \Gamma} \subset J_{\{n\}* \Sigma}$.
Since $\T (J_{\{n\}* \Sigma})=\{\Gin(J_{\{n\}* \Sigma})\}$
Lemma \ref{itumono2} (ii) says that
$J' \subset \Gin(J_{\{n\}* \Sigma})$.
Then Lemma \ref{cone} says that $\Gamma' \supset \{n\} * \dele \Sigma$.
We will show 
$$\Gamma' \supset \{n\} * \dele \Gamma.$$

Since $J_{\Gamma'}$ and $J_{\{n\} * \dele \Gamma}$ have the same Hilbert function,
the above inclusion implies $\Gamma' = \{n\} * \dele \Gamma$
and $|\T(J_{\{n\}* \Gamma})|=1$.
Also, to prove the above inclusion,
what we must prove is $\Gamma' \supset \{ \{n\} \cup S: S \in {\dele \Gamma }_k\}$
for $k=0,1,\dots,d-1$.

For $k=0,1,\dots,d-2$,
we have $\Gamma' \supset \{ \{n\} \cup S: S \in {\dele \Gamma }_k\}$ since $\Gamma_k=\Sigma_k$
and $\Gamma' \supset \{n\} * \dele \Sigma$.
Thus we will show the case $k=d-1$.
Since $\{ \{n\} *\Gamma \}_d= \{ \{n\} \cup S: S \in \Gamma_{d-1}\}$,
we have 
$$ (J_{\{n\} * \Gamma})_{d+1}= ( e_n J_{\Gamma\cup \{n\}}+(e_1,\dots,e_{n-1})^{d+1})_{d+1}.$$
Let $\varphi \in \GL$,
$\sigma$ a term order
and $e_t$ the monomial which does not belong to $\init_\sigma(\varphi(e_1,\dots,e_{n-1}))$.
Then,
for any $e_n \wedge f \in e_n J_{\Gamma \cup \{n\}}$,
we have $f \in (e_1,\dots,e_{n-1})^d$ and $e_k \wedge f \in (e_1,\dots,e_{n-1})^{d+1}$
for $k=1,2,\dots,n-1$.
Hence we have $e_k \wedge \varphi(f) \in \varphi( e_n J_{\Gamma\cup \{n\}}+(e_1,\dots,e_{n-1})^{d+1})_{d+1}$ for all $k \in [n]$.
This fact says that
\begin{eqnarray*}
&&\hspace{-50pt}\init_\sigma (\varphi( e_n J_{\Gamma\cup \{n\}}+(e_1,\dots,e_{n-1})^{d+1}))_{d+1}\\
&&\hspace{50pt}\supset
 ( e_t \init_\sigma(\varphi(J_{\Gamma\cup \{n\}}))+\init_\sigma (\varphi(e_1,\dots,e_{n-1}))^{d+1})_{d+1}.
\end{eqnarray*} 
Since $J'$ is strongly stable and since $\T(J_{\Gamma \cup \{n\}})=\{ \Gin(J_{\Gamma \cup \{n\}})\}$,
the above inclusion says that
$J'_{d+1} \supset  (e_n \Gin(J_{ \Gamma\cup \{n\}}) + (e_1,\dots,e_{n-1})^{d+1})_{d+1}$.
Thus we have
 $$\Gamma'_d \subset \{\{ n\} \cup S: S \in \Gamma_{d-1}\}.$$
However, since the cardinalities of both sides of the above inclusion are same,
we have $\Gamma' \supset \{ \{n\} \cup S: S \in {\dele \Gamma }_{d-1}\}$ as desired.
\end{proof}

\begin{lemma} \label{graphcone}
Let $G$ be a graph on $[n]$.
Assume that $G$ is a near cone with respect to $v \in [n]$.
If $|\T(J_{ \F (G-v) \cup \{v\}})|=1$,
then one has $|\T(J_{\F (G)})|=1$.
\end{lemma}

\begin{proof}
By Lemma \ref{itumono2} (iii),
we may assume that $v=n$.
We will show
\begin{eqnarray}
\{S \in \F(G): |S| \geq 3\} =\{S \in \{n\} * \F (G-n): |S| \geq 3\}. \label{new}
\end{eqnarray}
Note that the above equation and Lemma \ref{transcone} immediately imply $|\T(J_{\F (G)},d)|=1$ for $d \geq 3$.

The inclusion $\F(G) \subset\{n\} * \F (G-n)$ is true for an arbitrary graph $G$.
Let $S \in \{n\} * \F (G-n)$ with $|S| \geq 3$.
We will show $S \in \F(G)$.
If $n \not\in S$ then we have $S \in \F(G-n) \subset \F(G)$.
Otherwise, we have $S \setminus \{n\} \in \F(G-n) \subset \F (G)$.
Set $T= S \setminus \{n\}$.
Then we have $\{i,j\} \in G$ for any $\{i,j\} \subset T$.
Since $G$ is a near cone with respect to $n$ and $|T|>1$,
we have $\{i,j\} \in G$ for any $\{i,j\} \subset T \cup \{n\}$.
Thus $S=\{n\} \cup T \in \F (G)$.
Hence the equation (\ref{new}) follows.

It remains to prove $|\T(J_{\F (G)},2)|=1$.
Let $J' \in \T(J_{\F(G)})$ and $\Gamma'$ the simplicial complex with $J'=J_{\Gamma'}$.
We already show that $\F(G) \supset \{S \in \{n\} * \F(G-n):|S| =3\}$.
Then we have $\Gamma' \supset \dele {G-n}$ since $\T(J_{\{n\}*\F(G-n)})=\{J_{ \{n\}*\dele {\F(G-n)}}\}$
by the assumption and Lemma \ref{transcone}.
Let $H$ be the subgraph of $G$ with the edge set $\{ \{i,n\}: \{i,n\} \in G\}$.
Then $H$ is isomorphic to a shifted graph.
Thus we have $|\T (J_H)|=|\Gins(J_H)|=1$ and $\dele H$ have the edge set $\{\{n,n-j\}:j=1,2,\dots,\deg_G(n)\}$.
Since $H\subset \F(G)$,
Lemma \ref{itumono2} (ii) says that
$\dele H \subset \Gamma'$.
Hence we have
\begin{eqnarray}
\{S \in \Gamma': |S|=2\} \supset  \dele{G-n} \cup \{\{n,n-j\}:j=1,2,\dots,\deg_G(n)\}. \label{hosiQ}
\end{eqnarray}
Since $(J_{\F(G)})_2=(J_G)_2$,
the cardinality of the left-hand side of the above inclusion is equal to the number of edges in $G$.
On the other hand,
the number of edges in $G$ is equal to the sum of the number of edges in $G-n$ and $\deg_G(n)$.
Thus the inclusion (\ref{hosiQ}) is an equation.
Hence we have $|\T(J_{\F(G)},2)|=1$.
\end{proof}

Let $G$ be a graph on $[n]$.
Write $\overline{G}$ for the complementary graph of $G$.
By Corollary \ref{unique},
to prove $|\Gins(J_G)|=1$,
what we must prove is $\Delta^{\mathrm{lex}}(G)= \dele {G}$.
On the other hand,
as we saw in Example \ref{reisann},
since $<_{\mathrm{lex}^{-1}}$ is the reverse lexicographic order induced by $e_1<\cdots <e_n$,
$\Delta^{\mathrm{lex}}(G)$ can be computed from $\Delta^e(\overline G)$ by using Lemma \ref{complement}.
Now, we note the relation between $\Delta^{\mathrm{lex}}(G)$ and $\dele {\overline G}$.

Let $G$ be a graph on $[n]$ and
$f_1(G)$ the numbers of edges in $G$.
For an integer $k=1,2,\dots,n$,
define
$$\mathrm{max}_{ \geq k} (G) =|\{ \{i,j\}\in G:\max\{i,j\} \geq k\}|.$$
We also define $\max_{ \leq k}(G)$ and $\min_{\leq k}(G)$ in the same way.
Then, by a simple counting, we have
\begin{eqnarray}
\mathrm{max}_{\geq n+1-k}( \Delta^{\mathrm{lex}}(G))= {n \choose 2} - {n-k \choose 2}
-\Big\{ f_1(\overline G)-\mathrm{max}_{\leq n-k}( \overline{ \Delta^{\mathrm{lex}}(G)})\Big\}. \label{nidle}
\end{eqnarray}
Let $\pi : [n] \to [n]$ be the permutation defined by $\pi(j)=n+1-j$ for $j=1,2,\dots,n$.
Since $<_{\mathrm{lex}^{-1}}$ is the reverse lexicographic order induced by $e_1<\cdots <e_n$,
we have $\pi (\Delta^{\mathrm{lex}^{-1}}( \overline {G}))= \dele {\overline G}$ and 
\begin{eqnarray}
\mathrm{max}_{\leq n-k}( \Delta^{\mathrm{lex}^{-1}}( \overline G))= \mathrm{min}_{\geq k+1}( \dele {\overline G}). \label{spark}
\end{eqnarray}
Recall that Lemma \ref{complement} says that $\overline {\Delta^{\mathrm{lex}}(G)}= \Delta^{\mathrm{lex}^{-1}}(\overline G)$.
Thus equations (\ref{nidle}) and (\ref{spark}) yield the following relation.
\begin{eqnarray}
\mathrm{max}_{\geq n+1-k}( \Delta^{\mathrm{lex}}(G))= {n \choose 2} - {n-k \choose 2}
-\Big\{ f_1(\overline G)-\mathrm{min}_{\geq k+1}( \dele {\overline{G}})\Big\}. \label{ufo}
\end{eqnarray}

Let $a$ and $b$ be positive integers.
We write $K_{a,b}$ for the complete bipartite graph of size $a,b$
and $K_a \cup K_b$ for the disjoint union of two complete graphs $K_a$ of size $a$ and $K_b$ of size $b$,
where the vertex set of $K_{a,b}$ and $K_{a} \cup K_b$ is $[a+b]$.
Note that $K_a \cup K_b$ is the complementary graph of $K_{a,b}$.
The exterior algebraic shifted graph
$\dele {K_{a,b}}$ was computed by Kalai.

\begin{lemma}[{\cite[Theorem 6.1]{Kh}}] \label{bipargin}
Let $K_{a,b}$ be the complete bipartite graph of size $a,b$, where $a \leq b$.
Set $n=a+b$.
One has
\begin{eqnarray*}
\mathrm{max}_{\geq n+1- k}(\dele{K_{a,b}})=
\left \{
\begin{array}{l}
{n \choose 2} -{n-k \choose 2} -{k \choose 2}= kn-k^2,\ \mathrm{if} \ k\leq a, \\
ab,\ \ \  \ \ \ \ \hspace{107pt}        \mathrm{if} \ k>a .
\end{array}
\right.
\end{eqnarray*}
\end{lemma}

On the other hand,
we can compute $\dele {K_a \cup K_b}$ by using \cite[Theorem 1.1]{N}.

\begin{lemma} \label{disjgin}
Let $K_a \cup K_b$ be the disjoint union of two complete graphs $K_a$ and $K_b$, where $a \leq b$.
Set $n=a+b$.
Then one has
\begin{eqnarray}
\mathrm{min}_{\geq n+1- k}(\dele{K_a \cup K_b})=
\left \{
\begin{array}{l}
{k \choose 2},\  \ \ \hspace{87pt} \mathrm{if} \ k\leq b, \\
|f_1(K_a \cup K_b)| - {n-k \choose 2},\        \mathrm{if} \ k>b .
\end{array}
\right. \label{mike}
\end{eqnarray}
\end{lemma}

\begin{proof}
Since $K_a$ and $K_b$ are complete graphs,
we have $\dele {K_a} =K_a$ and $\dele {K_b}=K_b$ by Lemma \ref{itumono} (ii),
where we assume that $K_a$ is on $[a]$ and $K_b$ is on $[b]$.
Let 
$$ h_k= |\{ \{i,j\} \in \dele {K_a \cup K_b}: \max\{i,j\}= n+1-k\}|.$$
Then, it follows from \cite[Theorem 1.1]{N} that
\begin{eqnarray*}
h_k&=&|\{ \{i,j\} \in \dele {K_a}: \max\{i,j\}= a+1-k\}|\\
&&\hspace{70pt} +|\{ \{i,j\} \in \dele {K_b}: \max\{i,j\}= b+1-k\}|\\
&=&
\left\{
\begin{array}{ll}
(a-k) + (b-k), \ \ & \mbox{ if } k \leq a,\\
(b-k), &\mbox{ if } b \geq k > a,\\
0,   &\mbox{ if } k >b. 
\end{array}
\right.
\end{eqnarray*}
Since 
\begin{eqnarray*}
&&\mathrm{min}_{ \geq n+1 -k} ( \dele {K_a \cup K_b})\\
&&=\sum_{l=1}^{k-1} |\{ \{i,j\} \in \dele{K_a \cup K_b}: \max\{i,j\}=n+1-l,\ \min\{i,j\} \geq n+1-k\}|\\
&&=\sum_{l=1}^{k-1} \min\{k-l,h_l\},
\end{eqnarray*}
a routine computation yields (\ref{mike}).
\end{proof}

Let $G$ be a graph on $[n]$.
A vertex $v \in [n]$ is called an \textit{isolated vertex} of $G$ if $\deg_G(v)=0$.
Since $\Delta^\sigma(G)$ only depends on the combinatorial type of $G$ together with the characteristic of the field
(see \cite[\S 2]{K}),
$|\Gins(J_G)|$ does not change by deleting isolated vertices form $G$.

\begin{theorem}\label{extgraph}
If $G$ is a semi-complete bipartite graph
or a disjoint union of two semi-complete graphs,
then one has $|\Gins (J_G)|=1$.
\end{theorem}

\begin{proof}
We may assume that $G$ has no isolated vertices.
Moreover, by Lemma \ref{complement},
we have $|\Gins(J_G)|=1$ if and only if $|\Gins(J_{\overline G})|=1$.
Thus we may assume that $G$ is the complete bipartite graph $K_{a,b}$
for some positive integers $a$ and $b$.
Also, by Corollary \ref{unique},
what we must prove is $\Delta^{\mathrm{lex}}(K_{a,b})= \dele {K_{a,b}}$.
Now, by Lemmas \ref{bipargin} and \ref{disjgin} together with the equation (\ref{ufo}),
we have
$$
\mathrm{max}_{\geq n+1-k}( \Delta^{\mathrm{lex}}(K_{a,b}))=\mathrm{max}_{\geq n+1-k} ( \dele { K_{a,b}}) \ \ \mbox{ for } k=1,2,\dots,n.
$$
Hence we have $\max_{\leq k}(\Gin_{\mathrm{lex}}(J_{K_{a,b}}),2)=\max_{\leq k} (\Gin(J_{K_{a,b}}),2)$ for all $k$.
Thus we have $\Delta^{\mathrm{lex}}(K_{a,b})= \dele {K_{a,b}}$ by Lemma \ref{equivalent}.
\end{proof}

Let $I$ be a homogeneous ideal in the polynomial ring $R$.
The \textit{regularity} of $I$ is the integer
$\reg (I)=\max\{ d: \beta_{ii+d}(I) \ne 0 \mbox{ for some } i\}.$

\begin{lemma}\label{regularity}
Let $J$ be a monomial ideal in $E$ with $d= \mathrm{reg}( J^*)$ or a homogeneous ideal in $R$ with $\reg(J)=d$.
If $|\T(J,k)|=1$ for all $k \leq d$ then $|\T(J)|=1$.
\end{lemma}

\begin{proof}
We will show the case $J \subset E$.
(The proof for the case $J \subset R$ is same.)
Let $J' \in \T(J)$.
It follows from \cite[Corollary 3.6 and Theorem 7.1]{H} that
$\mathrm{reg}(J^*)$ is equal to the highest degree of
monomials belonging to the set of minimal monomial generators of $\Gin(J)$.
Then, for any $J' \in \T(J)$,
the assumption says that $J' \supset \Gin(J)$.
Since $J'$ and $\Gin (J)$ have the same Hilbert function,
we have $J'=\Gin(J)$.
Hence $\T(J)=\{\Gin(J)\}$.
\end{proof}

Finally, we require the next lemma.
A graph $G$ is said to be \textit{chordal}
if every induced cycle of $G$ has length $3$,
where an induced cycle of $G$ is a cycle of $G$ which is an induced subgraph of $G$.

\begin{lemma} [{Fr\"oesberg \cite{F}}] \label{froes}
A graph $G$ is chordal if and only if $\mathrm{reg}({J_{\F(G)}}^*) =2$.
\end{lemma}

\begin{theorem} \label{extflag}
Let $G$ be a $k$-near cone of a semi-complete bipartite graph
or of a disjoint union of two semi-complete graphs.
Then one has $|\T (J_{\F(G)})|=1$.
\end{theorem}

\begin{proof}
By Lemma \ref{graphcone},
we may assume that 
$G$ is either a semi-complete bipartite graph or a disjoint union of two semi-complete graphs.
Then Corollary \ref{unique} and Theorem \ref{extgraph} say $|\T (J_{\F(G)},2)|=1$.
If $G$ is a semi-complete bipartite graph,
then we have $|\T (J_{\F(G)})|= |\T (J_{G})|=1$.
Otherwise, $G$ is a chordal graph.
Since $\reg({J_{\F(G)}}^*)=2$  by Lemma \ref{froes} and $|\T(J_{\F(G)},2)|=1$,
we have $|\T (J_{\F(G)})|=1$
by  Lemma \ref{regularity}.
\end{proof}

\section{The proof of Theorem \ref{main1}}

In this section,
we will give a proof of Theorem \ref{main1}.
First, we will show (iv) $\Rightarrow$ (v).

\begin{proposition} \label{threegraphs}
Let $G$ be a graph on $[n]$.
If $G$ or $\overline G$ contains one of the graphs $(a)$, $(b)$ and $(c)$ given in \S 1 as an induced subgraph,
then one has $|\T (J_G)|>1$.
\end{proposition}

\begin{proof}
First, we will consider the graph (a).
Let $H_a$ be the graph with the edge set $\{\{1,2\},\{1,3\},\{3,4\}\}$.
Recall that Lemma \ref{complement} says that $|\Gins(J_G)|=1$ if and only if $|\Gins(J_{\overline G})|=1$ for an arbitrary graph $G$.
Then, by Corollary \ref{unique} and Lemma \ref{induced}, to prove the statement,
it suffices to show that $|\T (J_{\overline H_a})|>1$.
However, $J_{\overline H_a}$ is the ideal generated by $\{e_1 \wedge e_2,e_1 \wedge e_3, e_3 \wedge e_4\}$.
We already proved $|\T (J_{\overline H_a})|>1$ in Example \ref{rei}.
Thus we have $|\T(J_G)|>1$ if $G$ or $\overline G$ contains the graph (a) as an induced subgraph.
\medskip

Next, we will consider (b) and (c).
Let $H_b$ be the graph with the edge set $\{\{1,2\},\{3,4\},\{3,5\}\}$
and $H_c$ be the graph with the edge set $\{\{1,2\},\{3,4\},\{5,6\}\}$.
Consider graphs $H_b'$ and $H_c'$ defined by
$$J_{H_b'}= \init_{\mathrm{lex}} ( \varphi_{2,4} (J_{H_b})) \ \ \ \ \ 
\mbox{and} \ \ \ \ \ 
J_{H_c'}= \init_{\mathrm{lex}} ( \varphi_{3,6} (J_{H_c})),$$
where $\varphi_{i,j}$ is the matrix defined in Example \ref{rei}.
Then the graph $H_b'$ contains an induced subgraph which is isomorphic to the graph $H_a$
and  the graph $H_c'$ contains an induced subgraph which is isomorphic to the graph $H_b$.
Thus we have $|\T (J_{H_b},2)|>1$ and $|\T (J_{H_c},2)|>1$
by Lemma \ref{induced}.
Hence the claim follows from Lemmas \ref{complement} and \ref{induced}.
\end{proof}

Next,
we introduce some lemmas to prove (v) $\Rightarrow$ (vi) of Theorem \ref{main1}.
Let $G$ be a graph on $[n]$.
A \textit{proper connected component} of $G$ is a connected component of $G$
which is not an isolated vertex.

\begin{lemma} \label{various}
Let $G$ be a graph on $[n]$.
Assume that $G$ and $\overline G$ contain none of the graphs $(a)$, $(b)$ and $(c)$ as an induced subgraph.
Then
\begin{itemize}
\item[(i)] $G$ does not contain an induced cycle of length $\geq 5$;
\item[(ii)] if $G$ contains more than two proper connected components, then $G$ is a disjoint union of two semi-complete graphs;
\item[(iii)] if $G$ is a connected bipartite graph, then $G$ is a complete bipartite graph.
\end{itemize}
\end{lemma}

\begin{proof}
(i) Suppose that $G$ contains an induced cycle $\{a_1,a_2\},\{a_2,a_3\},\dots,\{a_t,a_1\}$ with $t \geq 5$.
Then the induced subgraph $G_{\{a_1,a_2,a_3,a_4\}}$ of $G$ on $\{a_1,a_2,a_3,a_4\}$ is isomorphic to the graph (a).
This contradicts the assumption.
\medskip

(ii) Suppose that $G$ contains more than two proper connected components and $G$ is not a disjoint union of two semi-complete graphs.
Then $G$ contains either (b) or (c) as an induces subgraph.
\medskip

(iii) 
Suppose that $G$ is not a complete bipartite graph.
Then there exist $i,j \in [n]$ such that
$\{i,j\} \not \in G$ and $G$ has a shortest path $i=a_1,\dots,a_{k-1},a_k=j$ from $i$ to $j$. 
Since $G$ is bipartite, the above path has length at least $3$,
and therefore we have $k \geq 4$.
Then $G_{\{a_1,\dots,a_4\}}$ is isomorphic to the graph (a).
\end{proof}

\begin{lemma} \label{twotype}
Let $G$ be a connected graph on the vertex set $A \cup B=[n]$ with $A \cap B = \emptyset$.
Assume  that $\{u,v\} \in G$ for all $u,v \in A$ and $ \{u',v'\} \not \in G$ for all $u',v' \in B$.
If $G$ and $\overline G$ contain none of the graphs $(a)$, $(b)$ and $(c)$ as an induced subgraph,
then there exists a vertex $a \in A$ such that $G$ is a near cone with respect to $a$.
\end{lemma}

To prove Lemma \ref{twotype},
we require the next lemma.

\begin{lemma}\label{totyuu}
With the same notation as in Lemma \ref{twotype},
assume that $V\! =\! \{v_1,\dots,v_k\}\\
 \subset A$ and $W=\{u_1,\dots,u_k\} \subset B$
are subsets which satisfy $\{u_i,v_j\} \in G$ if and only if $1 \leq i \leq j\leq k$.
If $G$ is not a near cone with respect to $v_k$,
then there exist $v_{k+1} \in A\setminus V$ and $u_{k+1}\in B\setminus W$
such that $(\mathrm{i})$ $\{u_t,v_{k+1}\} \in G$ for all $t=1,2,\dots,k+1$ and
$(\mathrm{ii})$ $\{u_{k+1},v_t\} \not \in G$ for all $t=1,2,\dots,k$.
\end{lemma}

\begin{proof}
Since $G$ is not a near cone with respect to $v_k$ and  $\{u,v\} \in G$ for all $u,v \in A$,
there exists a vertex $u_{k+1} \in B \setminus W$ such that $\{ u_{k+1},v_k\} \not \in G$.
First, we will show that this $u_{k+1}$ satisfies the condition (ii).
For each $t=1,2,\dots,k-1$, if $\{u_{k+1},v_t\} \in G$,
then the induced subgraph $G_{\{u_k,u_{k+1},v_t,v_k\}}$ is isomorphic to the graph (a).
This contradicts the assumption.
Thus we have $\{u_{k+1},v_t\} \not \in G$ for all $t=1,2,\dots,k$.

Next, since $G$ is connected and $ \{u',v'\} \not \in G$ for all $u',v' \in B$,
there exists a vertex $v_{k+1} \in A \setminus V$ such that $\{v_{k+1},u_{k+1}\} \in G$.
We will show that this $v_{k+1}$ satisfies the condition (i).
For each  $t=1,2,\dots,k$,
if $\{u_t,v_{k+1}\} \not\in G$ then the induced subgraph
$G_{\{u_t,u_{k+1},v_t,v_{k+1}\}}$ is isomorphic to the graph (a).
This contradicts the assumption.
Thus we have  $\{u_j,v_{k+1}\} \in G$ for all $t=1,2,\dots,k+1$.
\end{proof}

\begin{proof}[Proof of Lemma \ref{twotype}]
If $B= \emptyset$ then the statement is obvious.
Assume that $B \ne \emptyset$.
Let $u_1 \in B$.
Since $G$ is connected and $\{u,v\} \not \in G$ for all $u,v \in B$,
there exists a vertex $v_1 \in A$ such that $\{u_1,v_1\} \in G$.
Then subsets $V_1 =\{v_1\} \subset A$ and $W_1= \{u_1\} \subset B$ satisfy the conditions of 
Lemma \ref{totyuu}.
If $G$ is not a near cone w.r.t.\ $v_1$,
then Lemma \ref{totyuu} says that there exist $v_2 \in A\setminus V_1$ and $u_2 \in B \setminus W_1$
such that $V_2=\{v_1,v_2\}$ and $W_2=\{u_1,u_2\}$ satisfy the conditions of Lemma \ref{totyuu}.
Thus, arguing inductively,
there exists a vertex $a \in A$ such that $G$ is a near cone w.r.t.\ $a$.
\end{proof}

We already proved that
if a graph $G$ satisfies the condition of Lemma \ref{various},
then the length of every induced cycle of $G$ is either 3 or 4.
We will consider the case that $G$ contains an induced cycle of length $4$.

\begin{lemma} \label{length4}
With the same notation as in Lemma \ref{various}.
If $G$ is connected and contains an induced cycle of length $4$,
then $G$ is a bipartite graph or
there exists a vertex $v_0 \in [n]$ such that $G$ is a near cone with respect to $v_0$.
\end{lemma}

\begin{proof}
We assume that $G$ contains an induced cycle of length $4$ on $A=\{a_1,\dots,a_4\}\subset [n]$.
Then, for any vertex $v \in [n]\setminus A$,
the induced subgraph $G_{\{v\} \cup A}$ must be isomorphic to one of the following graphs.

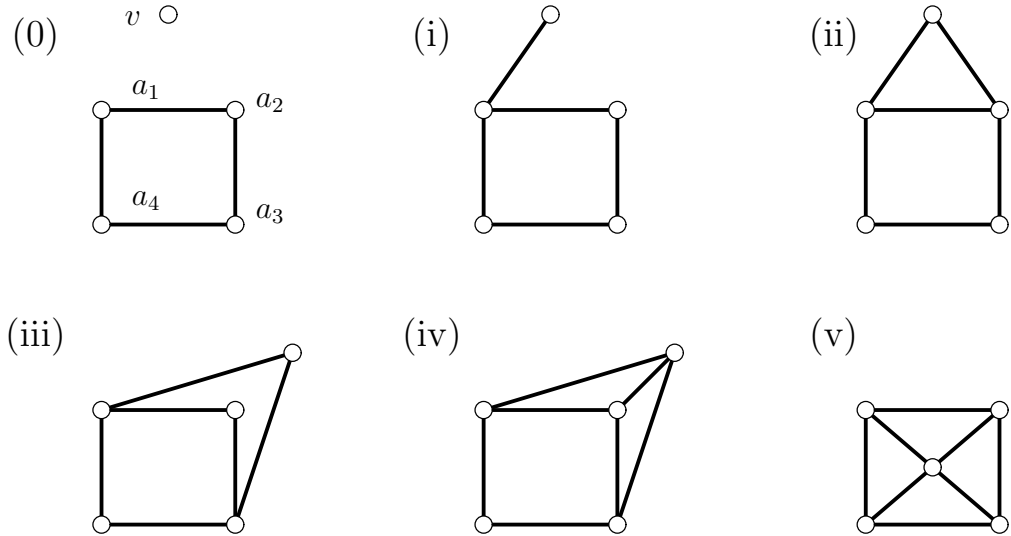
\begin{picture}(36.47, 22.63)(12.25, -30.63)
\put (11, -30){\makebox (0, 0){\large{(0)}}}%
\put (48, -23){\makebox (0, 0){$v$}}%
\put (53, -50){\makebox (0, 0){$a_1$}}%
\put (100, -55){\makebox (0, 0){$a_2$}}%
\put (53, -92){\makebox (0, 0){$a_4$}}%
\put (100, -97){\makebox (0, 0){$a_3$}}%

\put (161, -30){\makebox (0, 0){\large{(i)}}}%
\put (313, -30){\makebox (0, 0){\large{(ii)}}}%

\special{pn 20}%
\special{pa 500 800}%
\special{pa 1200 800}%
\special{fp}%

\special{pn 20}%
\special{pa 500 1400}%
\special{pa 1200 1400}%
\special{fp}%

\special{pn 20}%
\special{pa 1200 800}%
\special{pa 1200 1400}%
\special{fp}%

\special{pn 20}%
\special{pa 500 800}%
\special{pa 500 1400}%
\special{fp}%

\special{pn 8}%
\special{sh 0.000}%
\special{ar 850 300 46 46  0.0000000 6.2831853}%

\special{pn 8}%
\special{sh 0.000}%
\special{ar 500 800 46 46  0.0000000 6.2831853}%

\special{pn 8}%
\special{sh 0.000}%
\special{ar 500 1400 46 46  0.0000000 6.2831853}%

\special{pn 8}%
\special{sh 0.000}%
\special{ar 1200 800 46 46  0.0000000 6.2831853}%

\special{pn 8}%
\special{sh 0.000}%
\special{ar 1200 1400 46 46  0.0000000 6.2831853}%

\special{pn 20}%
\special{pa 2500 800}%
\special{pa 3200 800}%
\special{fp}%

\special{pn 20}%
\special{pa 2500 1400}%
\special{pa 3200 1400}%
\special{fp}%

\special{pn 20}%
\special{pa 3200 800}%
\special{pa 3200 1400}%
\special{fp}%

\special{pn 20}%
\special{pa 2500 800}%
\special{pa 2500 1400}%
\special{fp}%

\special{pn 20}%
\special{pa 2850 300}%
\special{pa 2500 800}%
\special{fp}%

\special{pn 8}%
\special{sh 0.000}%
\special{ar 2850 300 46 46  0.0000000 6.2831853}%

\special{pn 8}%
\special{sh 0.000}%
\special{ar 2500 800 46 46  0.0000000 6.2831853}%

\special{pn 8}%
\special{sh 0.000}%
\special{ar 2500 1400 46 46  0.0000000 6.2831853}%

\special{pn 8}%
\special{sh 0.000}%
\special{ar 3200 800 46 46  0.0000000 6.2831853}%

\special{pn 8}%
\special{sh 0.000}%
\special{ar 3200 1400 46 46  0.0000000 6.2831853}%

\special{pn 20}%
\special{pa 4500 800}%
\special{pa 5200 800}%
\special{fp}%

\special{pn 20}%
\special{pa 4500 1400}%
\special{pa 5200 1400}%
\special{fp}%

\special{pn 20}%
\special{pa 5200 800}%
\special{pa 5200 1400}%
\special{fp}%

\special{pn 20}%
\special{pa 4500 800}%
\special{pa 4500 1400}%
\special{fp}%

\special{pn 20}%
\special{pa 4850 300}%
\special{pa 4500 800}%
\special{fp}%

\special{pn 20}%
\special{pa 4850 300}%
\special{pa 5200 800}%
\special{fp}%


\special{pn 8}%
\special{sh 0.000}%
\special{ar 4850 300 46 46  0.0000000 6.2831853}%

\special{pn 8}%
\special{sh 0.000}%
\special{ar 4500 800 46 46  0.0000000 6.2831853}%

\special{pn 8}%
\special{sh 0.000}%
\special{ar 4500 1400 46 46 0.0000000 6.2831853}%

\special{pn 8}%
\special{sh 0.000}%
\special{ar 5200 800 46 46  0.0000000 6.2831853}%

\special{pn 8}%
\special{sh 0.000}%
\special{ar 5200 1400 46 46  0.0000000 6.2831853}%

\end{picture}%
\vspace{90pt}

\begin{picture}(36.47, 22.63)(12.25, -30.63)
\put (11, -30){\makebox (0, 0){\large{(iii)}}}%

\put (161, -30){\makebox (0, 0){\large{(iv)}}}%
\put (313, -30){\makebox (0, 0){\large{(v)}}}%

\special{pn 20}%
\special{pa 500 800}%
\special{pa 1200 800}%
\special{fp}%

\special{pn 20}%
\special{pa 500 1400}%
\special{pa 1200 1400}%
\special{fp}%

\special{pn 20}%
\special{pa 1200 800}%
\special{pa 1200 1400}%
\special{fp}%

\special{pn 20}%
\special{pa 500 800}%
\special{pa 500 1400}%
\special{fp}%

\special{pn 20}%
\special{pa 1500 500}%
\special{pa 500 800}%
\special{fp}%

\special{pn 20}%
\special{pa 1500 500}%
\special{pa 1200 1400}%
\special{fp}%

\special{pn 8}%
\special{sh 0.000}%
\special{ar 1500 500 46 46  0.0000000 6.2831853}%

\special{pn 8}%
\special{sh 0.000}%
\special{ar 500 800 46 46  0.0000000 6.2831853}%

\special{pn 8}%
\special{sh 0.000}%
\special{ar 500 1400 46 46  0.0000000 6.2831853}%

\special{pn 8}%
\special{sh 0.000}%
\special{ar 1200 800 46 46  0.0000000 6.2831853}%

\special{pn 8}%
\special{sh 0.000}%
\special{ar 1200 1400 46 46  0.0000000 6.2831853}%

\special{pn 20}%
\special{pa 2500 800}%
\special{pa 3200 800}%
\special{fp}%

\special{pn 20}%
\special{pa 2500 1400}%
\special{pa 3200 1400}%
\special{fp}%

\special{pn 20}%
\special{pa 3200 800}%
\special{pa 3200 1400}%
\special{fp}%

\special{pn 20}%
\special{pa 2500 800}%
\special{pa 2500 1400}%
\special{fp}%

\special{pn 20}%
\special{pa 3500 500}%
\special{pa 2500 800}%
\special{fp}%

\special{pn 20}%
\special{pa 3500 500}%
\special{pa 3200 1400}%
\special{fp}%

\special{pn 20}%
\special{pa 3500 500}%
\special{pa 3200 800}%
\special{fp}%

\special{pn 8}%
\special{sh 0.000}%
\special{ar 3500 500 46 46  0.0000000 6.2831853}%

\special{pn 8}%
\special{sh 0.000}%
\special{ar 2500 800 46 46  0.0000000 6.2831853}%

\special{pn 8}%
\special{sh 0.000}%
\special{ar 2500 1400 46 46  0.0000000 6.2831853}%

\special{pn 8}%
\special{sh 0.000}%
\special{ar 3200 800 46 46  0.0000000 6.2831853}%

\special{pn 8}%
\special{sh 0.000}%
\special{ar 3200 1400 46 46 0.0000000 6.2831853}%

\special{pn 20}%
\special{pa 4500 800}%
\special{pa 5200 800}%
\special{fp}%

\special{pn 20}%
\special{pa 4500 1400}%
\special{pa 5200 1400}%
\special{fp}%

\special{pn 20}%
\special{pa 5200 800}%
\special{pa 5200 1400}%
\special{fp}%

\special{pn 20}%
\special{pa 4500 800}%
\special{pa 4500 1400}%
\special{fp}%

\special{pn 20}%
\special{pa 4850 1100}%
\special{pa 4500 800}%
\special{fp}%

\special{pn 20}%
\special{pa 4850 1100}%
\special{pa 5200 800}%
\special{fp}%

\special{pn 20}%
\special{pa 4850 1100}%
\special{pa 4500 1400}%
\special{fp}%

\special{pn 20}%
\special{pa 4850 1100}%
\special{pa 5200 1400}%
\special{fp}%


\special{pn 8}%
\special{sh 0.000}%
\special{ar 4850 1100 46 46  0.0000000 6.2831853}%

\special{pn 8}%
\special{sh 0.000}%
\special{ar 4500 800 46 46  0.0000000 6.2831853}%

\special{pn 8}%
\special{sh 0.000}%
\special{ar 4500 1400 46 46  0.0000000 6.2831853}%

\special{pn 8}%
\special{sh 0.000}%
\special{ar 5200 800 46 46  0.0000000 6.2831853}%

\special{pn 8}%
\special{sh 0.000}%
\special{ar 5200 1400 46 46  0.0000000 6.2831853}%
\end{picture}%
\vspace{90pt}

However, (i) and (ii) contains the graph (a) as an induced subgraph
and the complementary graph of (iv) is isomorphic to the graph (b).
Hence the induced subgraph $G_{\{v\} \cup A}$
must be isomorphic to one of the graphs (0), (iii) and (v).

Let 
\begin{eqnarray*}
X_0&=&\{ v \in [n] \setminus A: \{v,a_k\} \not \in G \mbox{ for all } a_k \in A\},\\
X_1&=&\{ v \in [n] \setminus A: \{v,a_k\} \in G \mbox{ for $k=1,3$ and } \{v,a_k\}\not \in G \mbox{ for $k=2,4$}\},\\
X_2&=&\{ v \in [n] \setminus A: \{v,a_k\} \not \in G \mbox{ for $k=1,3$ and } \{v,a_k\} \in G \mbox{ for $k=2,4$}\},\\
X_4&=&\{ v \in [n] \setminus A: \{v,a_k\} \in G \mbox{ for all } a_k \in A\}.
\end{eqnarray*}
Note that $X_0\cup X_1 \cup X_2 \cup X_4 = [n] \setminus A$.
The next claim easily follows.
\bigskip

\noindent
\textbf{[Claim]}
\begin{itemize}
\item[(I)] If $u,v \in X_4$, then $\{u,v \} \in G$.
\item[(II)] If $u \in X_1 \cup X_2$ and $v \in X_4$, then we have $\{u,v\} \in G$.
\item[(III)] If $u \in X_1\cup X_2$ and $v \in X_0$, then we have $\{u,v\} \not\in G$.
\item[(IV)] If  $u,v \in X_0$ then we have $\{u,v\} \not \in G$.
\item[(V)] For any $u \in X_0$, there exist $v \in X_4$ such that $\{u,v\} \in G$.
\end{itemize}

\begin{proof}[Proof of Claim]
We will show that if $G$ does not satisfy one of (I), (II), (III) and (IV)
then $G$ contains one of the graphs (a), (b) and (c) as an induced subgraph.
(I) If $\{u,v\} \not \in G$, then $\overline G_{ \{u,v\} \cup A}$ is isomorphic to the graph (c).
(II) We may assume that $u \in X_1$.
If $\{u,v\} \not \in G$, then $\overline G_{\{a_1,a_2,a_3,u,v\}}$ is isomorphic to
the graph (b).
(III) Assume that $u \in X_1$.
If $\{u,v\} \in G$, then  $G_{\{a_1,a_2,u,v\}}$ is isomorphic to the graph (a).
(IV) If $\{u,v\} \in G$, then $G_{\{a_1,a_2,a_3,u,v\}}$ is isomorphic to the graph (b).

Finally, we will show (V).
Let $u \in X_0$.
Since $G$ is connected, there exists a vertex $v \in [n] \setminus A$ such that $\{u,v\}\in G$.
Then (III) and (IV) say that $v \in X_4$. 
\end{proof}

Now, we return the proof of Lemma \ref{length4}.\medskip

\textbf{[Case 1]} We will show that if $X_4 \ne \emptyset$ then
there exists a vertex $v_0 \in X_4$ such that $G$ is a near cone with respect to $v_0$.

Consider the graph $H= G_{X_0 \cup X_4}$.
Then, [Claim] (I) and (V) say that $H$ is connected.
Furthermore, [Claim] (I) and (IV) say that $H$ satisfies the conditions of Lemma \ref{twotype}.
Thus there exists a vertex $v_0 \in X_4$ such that $H$ is a near cone with respect to $v_0$.
Then [Claim]  (II) says hat $G$ is a near cone with respect to $v_0$.
\medskip

\textbf{[Case 2]}
We will show that if $X_4 = \emptyset$ then $G$ is bipartite.

If $X_4 = \emptyset$ then $X_0$ is also empty by [Claim] (V).
Let $B=\{a_1,a_3\} \cup X_2$ and $C=\{a_2,a_4\} \cup X_1$.
We will show that $G$ is a bipartite graph with the bipartition $\{B,C\}$.
By the construction of $X_1$ and $X_2$, what we must prove is that
$\{u,v\}  \not \in G$ if $\{u,v\}\subset X_1$ or $\{u,v\} \subset X_2$.
We may assume that $\{u,v\} \subset X_1$.
Suppose that $\{u,v\} \in G$.
Then
$\overline G_{\{a_1,a_2,a_3,u,v\}}$ is isomorphic to the graph (b).
This contradicts the assumption.
Thus we have  $\{u,v\}  \not \in G$ if $\{u,v\} \subset  X_1$ or $\{u,v\} \subset X_2$.
Hence $G$ is a bipartite graph if $X_4 = \emptyset$.
\end{proof}

Next, we consider the case that $G$ does not contain
an induced cycle of length $4$.

\begin{lemma} \label{chordal}
With the same notation as in Lemma \ref{various}.
If $G$ is a connected chordal graph,
then there exists a vertex $v_0 \in [n]$ such that  $G$ is a near cone with respect to $v_0$.
\end{lemma}

\begin{proof}
If $G$ is a complete graph, then $G$ is a near cone w.r.t.\ any $v \in [n]$.
We assume that $G$ is not a complete graph.
Then there exist $i,j \in [n]$ such that $\{i,j\} \not \in G$.
Let 
\begin{eqnarray*}
A&=&\{v \in [n ]\setminus \{i,j\}: \{i,v\} \in G \mbox{ and } \{j,v\}\in G \}\\ 
B&=&\{v \in [n] \setminus \{i,j\}: \{u,v\} \not \in G \mbox{ for all } u \in [n] \setminus A\}\\
\mbox{ and }\ \ C&=&[n] \setminus (\{i,j\}\cup A \cup B).
\end{eqnarray*}
First, we will show that $A \ne \emptyset$ and $\{u,v \} \in G$ for all $u,v \in A$.

If $A$ is empty, then the shortest path from $i$ to $j$ has at least length $3$.
Then $G$ contains the graph (a) as an induced subgraph in the same way as
the proof of Lemma \ref{various} (iii).
Thus $A \ne \emptyset$.
On the other hand, if $u ,v \in A$,
then $\{i,u\},\{u,j\},\{j,v\},\{v,i\}$ is a cycle of $G$.
Since $\{i,j\} \not \in G$ and $G$ is chordal,
we have $\{u,v\} \in G$.
\medskip

Second, we will show that $\{u,v \} \in G$ for any $u \in A$ and all $v \in C$.

Let $u \in A$ and $v \in C$.
Suppose that $\{u,v\} \not \in G$.
Since $G$ is connected,
by the construction of $C$,
there exists a vertex $v' \in A\cup C \cup \{i,j\}$  such that $\{v,v'\} \in G$.
Since $v \not \in B$, we may assume that $v' \not \in A$.
Also, since $v \not \in A$, we have either $\{v,i\}\not \in G$ or $\{v,j\} \not \in G$.
If $\{v,i\}\in G$ or $\{v,j\}\in G$ then $G_{\{u,v,i,j\}}$ is isomorphic to the graph (a).
Thus we may assume that $\{v,i\}\not\in G$, $\{v,j\}\not\in G$
and $v'\in C$.
Then, since $v' \not \in A$,
we have either $\{v',i\}\not \in G$ or $\{v',j\} \not \in G$.
Hence the induced subgraph $G_{\{u,v,v',i,j\}}$ is isomorphic to one of the following graphs.
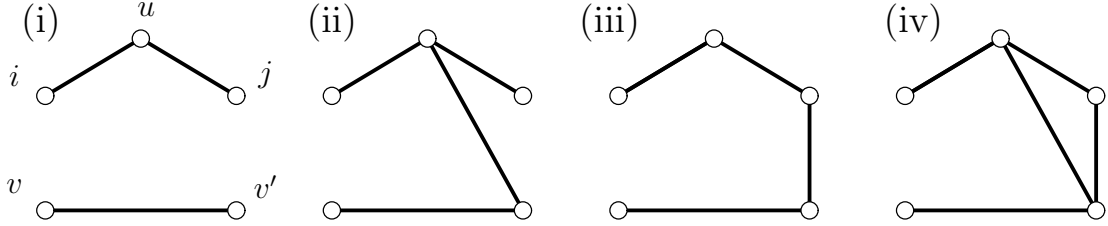
\vspace{-3pt}

\begin{picture}(36.47, 22.63)(12.25, -30.63)
\put (20, -30){\makebox (0, 0){\large{(i)}}}%
\put (60, -25){\makebox (0, 0){$u$}}%
\put (10, -50){\makebox (0, 0){$i$}}%
\put (105, -50){\makebox (0, 0){$j$}}%
\put (10, -92){\makebox (0, 0){$v$}}%
\put (105, -92){\makebox (0, 0){$v'$}}%

\put (130, -30){\makebox (0, 0){\large{(ii)}}}%
\put (235, -30){\makebox (0, 0){\large{(iii)}}}%
\put (350, -30){\makebox (0, 0){\large{(iv)}}}%

\special{pn 20}%
\special{pa 300 1400}%
\special{pa 1300 1400}%
\special{fp}%

\special{pn 20}%
\special{pa 800 500}%
\special{pa 1300 800}%
\special{fp}%

\special{pn 20}%
\special{pa 800 500}%
\special{pa 300 800}%
\special{fp}%

\special{pn 8}%
\special{sh 0.000}%
\special{ar 800 500 46 46  0.0000000 6.2831853}%

\special{pn 8}%
\special{sh 0.000}%
\special{ar 300 800 46 46  0.0000000 6.2831853}%

\special{pn 8}%
\special{sh 0.000}%
\special{ar 300 1400 46 46  0.0000000 6.2831853}%

\special{pn 8}%
\special{sh 0.000}%
\special{ar 1300 800 46 46  0.0000000 6.2831853}%

\special{pn 8}%
\special{sh 0.000}%
\special{ar 1300 1400 46 46  0.0000000 6.2831853}%

\special{pn 20}%
\special{pa 1800 1400}%
\special{pa 2800 1400}%
\special{fp}%

\special{pn 20}%
\special{pa 2300 500}%
\special{pa 2800 1400}%
\special{fp}%

\special{pn 20}%
\special{pa 2300 500}%
\special{pa 2800 800}%
\special{fp}%

\special{pn 20}%
\special{pa 2300 500}%
\special{pa 1800 800}%
\special{fp}%

\special{pn 8}%
\special{sh 0.000}%
\special{ar 2300 500 46 46  0.0000000 6.2831853}%

\special{pn 8}%
\special{sh 0.000}%
\special{ar 1800 800 46 46  0.0000000 6.2831853}%

\special{pn 8}%
\special{sh 0.000}%
\special{ar 1800 1400 46 46  0.0000000 6.2831853}%

\special{pn 8}%
\special{sh 0.000}%
\special{ar 2800 800 46 46  0.0000000 6.2831853}%

\special{pn 8}%
\special{sh 0.000}%
\special{ar 2800 1400 46 46  0.0000000 6.2831853}%

\special{pn 20}%
\special{pa 3300 1400}%
\special{pa 4300 1400}%
\special{fp}%

\special{pn 20}%
\special{pa 3800 500}%
\special{pa 4300 800}%
\special{fp}%

\special{pn 20}%
\special{pa 4300 800}%
\special{pa 4300 1400}%
\special{fp}%

\special{pn 20}%
\special{pa 3800 500}%
\special{pa 3300 800}%
\special{fp}%

\special{pn 20}%
\special{pa 3800 500}%
\special{pa 3300 800}%
\special{fp}%

\special{pn 8}%
\special{sh 0.000}%
\special{ar 3800 500 46 46  0.0000000 6.2831853}%

\special{pn 8}%
\special{sh 0.000}%
\special{ar 3300 800 46 46  0.0000000 6.2831853}%

\special{pn 8}%
\special{sh 0.000}%
\special{ar 3300 1400 46 46  0.0000000 6.2831853}%

\special{pn 8}%
\special{sh 0.000}%
\special{ar 4300 800 46 46  0.0000000 6.2831853}%

\special{pn 8}%
\special{sh 0.000}%
\special{ar 4300 1400 46 46  0.0000000 6.2831853}%

\special{pn 20}%
\special{pa 4800 1400}%
\special{pa 5800 1400}%
\special{fp}%

\special{pn 20}%
\special{pa 5300 500}%
\special{pa 5800 800}%
\special{fp}%

\special{pn 20}%
\special{pa 5800 800}%
\special{pa 5800 1400}%
\special{fp}%

\special{pn 20}%
\special{pa 5300 500}%
\special{pa 5800 1400}%
\special{fp}%

\special{pn 20}%
\special{pa 5300 500}%
\special{pa 4800 800}%
\special{fp}%

\special{pn 20}%
\special{pa 5300 500}%
\special{pa 4800 800}%
\special{fp}%

\special{pn 8}%
\special{sh 0.000}%
\special{ar 5300 500 46 46  0.0000000 6.2831853}%

\special{pn 8}%
\special{sh 0.000}%
\special{ar 4800 800 46 46  0.0000000 6.2831853}%

\special{pn 8}%
\special{sh 0.000}%
\special{ar 4800 1400 46 46  0.0000000 6.2831853}%

\special{pn 8}%
\special{sh 0.000}%
\special{ar 5800 800 46 46  0.0000000 6.2831853}%

\special{pn 8}%
\special{sh 0.000}%
\special{ar 5800 1400 46 46  0.0000000 6.2831853}%

\end{picture}
\vspace{90pt}

Then (i) is isomorphic to the graph (b)
and (ii), (iii) and (iv) contains the graph (a) as an induced subgraph.
This contradicts the assumption that $G$ does not contain the graphs (a) and (b) as an induced subgraph.
Thus we have $\{u,v\} \in G$ for any $u \in A$ and $v\in C$.
\medskip

Now, we will prove the statement.
Let $H= G_{A \cup B}$.
Since $G$ is connected, for any $b\in B$,
there exists a vertex $v_b \in A$ such that $\{b,v_b\} \in G$.
Also, we proved that $\{u,v\} \in G$ for all $u,v \in A$.
These facts say that $H$ is connected and
satisfies the conditions of Lemma \ref{twotype}.
Thus there exists a vertex $v_0 \in A$ such that $H$ is a near cone w.r.t.\ $v_0.$
We already proved $\{u,v\} \in G$ for any $u \in A$ and $v \in C$.
Since $v_0 \in A$, it follows that the graph $G$ is a near cone with respect to $v_0$, as desired.
\end{proof}

Now, we are in the position to give a proof of Theorem \ref{main1}.

\begin{proof}[Proof of Theorem \ref{main1}]
(i) $\Rightarrow$ (ii) is obvious and
(ii) $\Rightarrow$ (iii) follows from the facts that
$(J_{\F(G)})_2=(J_G)_2$ and $|\Gins(J_G,d)|=1$ for all $d \geq 3$.
Also,
(iii) $\Rightarrow$ (iv) follows from Corollary \ref{unique},
(iv) $\Rightarrow$ (v) is Proposition \ref{threegraphs}
and (vi) $\Rightarrow$ (i) is Theorem \ref{extflag}.

We will show (v) $\Rightarrow$ (vi).
Assume that $G$ and $\overline G$ contain none of the graphs $(a)$, $(b)$ and $(c)$ as an induced subgraph.
Let $G'$ be an induced subgraph of $G$ such that
$G$ is the $k$-near cone of a graph $G'$ and $G'$ is not a near cone w.r.t.\ $v$ for any $v \in [n]$.
If $G'$ contains no edges,
then $G$ is a $(k-1)$-near cone of a star modulo isolated vertices,
that is,
$G$ is a $(k-1)$-near cone of a union of the complete bipartite graph of size $1,t$
for some integer $t>0$ and isolated vertices.

Next, assume that $G'$ contains an edge.
By Lemma \ref{various},
if $G'$ has more than two proper connected components,
then $G$ satisfies the condition (vi) of Theorem \ref{main1}.
On the other hand, if $G'$ has only one proper connected component $H$,
then, by Lemma \ref{various}, $H$ does not have an induced cycle of length $\geq 5$.
Also, for any $v \in [n]$, since 
$G'$ is not a near cone w.r.t.\ $v$,
$H$ is not a near cone w.r.t.\ $v$.
Then Lemma \ref{chordal} says that $H$ is not a chordal graph.
Thus $H$ contains an induced cycle of length $4$.
Then Lemmas \ref{various} and \ref{length4} say that $H$ is a complete bipartite graph.
\end{proof}

\section{Edge ideals}

Let $K$ be a field of characteristic $0$
and $R=K[x_1,\dots,x_n]$ the polynomial ring in $n$ variables with each $\deg (x_i)=1$.
In this section, the proof of Theorem \ref{main2} will be given.
We split Theorem \ref{main2} into Theorem \ref{suff2} and Theorem \ref{suff}.

First, we will
prove the ``if" part of Theorem \ref{main2}.
Let $\Gamma$ be a simplicial complex on $[n]$.
The \textit{Stanley--Reisner} ideal $I_\Gamma \subset R$ of $\Gamma$ is the monomial ideal generated by 
all squarefree monomials $x_{i_1}\cdots x_{i_k}$ with $\{i_1,\dots,i_k\} \not \in \Gamma$. 
Let $I \subset R$ be a homogeneous ideal.
We say that $I$ has a \textit{linear resolution} if $I$ is generated in degree $d$
and $\mathrm{reg}(I)=d$.
For example,
by Lemma \ref{froes},
if $G$ is a chordal graph then $I_{\F(G)}=I(\overline G)$ has a linear resolution.

\begin{lemma}[{\cite[Theorem 1.2 and Corollary 2.2]{AHH}}] \label{linear}
Let $\Gamma$ be a simplicial complex.
If $I_{\Gamma}$ has a linear resolution,
then 
$$\beta_{ij}(I_\Gamma)= \beta_{ij}(\gin(I_{\Gamma})) = \beta_{ij} (  {\Gin (J_\Gamma)}^*)\ \ \ \ \mbox{ for all }i \mbox{ and }j.$$ 
\end{lemma}

Let $\mathcal{S}_2 \subset R_2$ be  the set of monomials in $R$ of degree $2$
and $\mathcal{M}_2 \subset E_2$ the set of monomials in $E$ of degree $2$.
Let $V$ be a $K$-vector space with basis $e_1,\dots,e_n$ and
$\varphi=(a_{ij}) \in GL_n(K)$.
For $k=1,2,\dots,n$,
define the map $\rho_{\varphi,k}:\mathcal{M}_2 \to \bigoplus_{i=1}^k V$ by
\[
\rho_{\varphi,k} (e_i\wedge e_j)=(a_{1j}e_i-a_{1i}e_j,\dots,a_{kj}e_i-a_{ki}e_j) \in \bigoplus_{i=1}^k V
\]
and $\phi_{\varphi,k}: \mathcal{S}_2 \to \bigoplus_{i=1}^k V$ by
$$\phi_{\varphi,k} (x_ix_j) = (a_{1j}e_i+a_{1i}e_j,\dots,a_{kj}e_i+a_{ki}e_j)\in \bigoplus_{i=1}^k V.$$

The next property is known.
(See \cite[Lemma 7.1]{K} and \cite[Lemma 1.8]{Mf}.)

\begin{lemma} \label{hyper}
Let $W_1 \subset E_2$ be a $K$-vector space spanned by monomials in $E$
and $W_2 \subset  R_2$ a $K$-vector space spanned by monomials in $R$.
Then, for a generic matrix $\varphi \in \GL$ and for $k=1,2,\dots,n$, one has
$$|\{ e_i \wedge e_j \not \in \Gin(W_1): \max\{i,j\} \geq k\}| =
\dim_K \mathrm{span} \{\rho_{\varphi,n+1-k}(e_i \wedge e_j):e_i \wedge e_j \not \in W_1 \}$$
and
$$|\{ x_i x_j \not \in \Gin(W_2): \max\{i,j\} \geq k\}| =
\dim_K \mathrm{span} \{\phi_{\varphi,n+1-k}(x_i x_j):x_i x_j \not \in W_2 \}.$$
\end{lemma}

For a graph $G$ on $[n]$,
set $\rho_{\varphi,k}(G)=\mathrm{span}\{\rho_{\varphi,k} (e_i \wedge e_j): \{i,j\} \in G\}$
and $\phi_{\varphi,k}(G)=\mathrm{span}\{\phi_{\varphi,k} (x_ix_j): \{i,j\} \in G\}$.
Then, by using the fact that $<_{\mathrm{lex}}^{-1}$ is the reverse lexicographic order
induced by $1<2< \dots <n$,
Lemmas \ref{symcomplement} and \ref{complement} together with Lemma \ref{hyper} say that,
for a generic matrix $\varphi \in \GL$, one has
\begin{eqnarray}
\mathrm{min}_{\leq k} ( \Gin_{\mathrm{lex}}(J_{\overline G}),2) = \dim_K \rho_{\varphi,k}(G) \ \ \ \ \mbox{for all } k \label{MG1}
\end{eqnarray}
and
\begin{eqnarray}
\mathrm{min}_{\leq k} ( \gin_{\mathrm{lex}}(I(G),2) = \dim_K \phi_{\varphi,k}(G) \ \ \ \ \mbox{for all } k. \label{MG2}
\end{eqnarray}
By using the above formula,
we can prove the next property.

\begin{lemma} \label{finite}
If $G$ is a bipartite graph,
then one has 
$$\mathrm{min}_{\leq k} ( \Gin_{\mathrm{lex}}(J_{\overline G}),2) = \mathrm{min}_{\leq k} ( \gin_{\mathrm{lex}}(I(G)),2) \ \ \ \mbox{ for all }k.$$
\end{lemma}

\begin{proof} (sketch)
The proof of this lemma essentially appeared in \cite[Lemma 2.2]{Mf}.
Thus we sketch the proof.
Assume that $G$ is a bipartite graph with the bipartition $\{A,B\}$.
Define the automorphism $\Phi:V \to V$ by
$\Phi(e_i)=e_i$ for $i \in A$, and $\Phi(e_i)=-e_i$ for $i \in B$.
Let $\Phi^{(k)}:\bigoplus_{i=1}^k V \to \bigoplus_{i=1}^k V$
be the automorphism
defined by $\Phi(u_1,\dots,u_k)=(\Phi(u_1),\dots,\Phi(u_k))$.
Then, for any $\varphi \in \GL$ and for any integer $k$,
one has $\Phi ^{(k)} (\rho_{\varphi,k}(G))= \phi_{\varphi,k}(G)$
since $G$ is bipartite.
Thus the claim follows from the equations (\ref{MG1}) and (\ref{MG2}).
\end{proof}

For any monomial $x_{i_1}x_{i_2}\cdots x_{i_k} \in R$ with $i_1 \leq i_2 \leq \cdots \leq i_k$,
define 
$$\alpha(x_{i_1}x_{i_2}\cdots x_{i_k})= x_{i_1}x_{i_2+1}\cdots x_{i_k+k-1}.$$
Also, for any strongly stable ideal $I\subset R$,
we write $\alpha(I)$ for the ideal generated by $\{ \alpha(u): u \in G(I)\}$,
where $G(I)$ is the set of minimal monomial generators of $I$ and we  assume that $\alpha(u) \in K[x_1,\dots,x_n]$ for all $u \in G(I)$.
Then, it follows from \cite[Lemmas 8.17 and 8.18]{H} that
\begin{eqnarray}
\beta_{ij}(I)=\beta_{ij}(\alpha(I)) \ \ \ \mbox{ for all } i \mbox{ and }j \label{SRM1}
\end{eqnarray}
and $\alpha(I)$ is a squarefree strongly stable ideal,
that is, there exists a strongly stable ideal $J \subset E$ such that  $ J^*=\alpha (I)$.
In particular,
if $\Gamma$ is a simplicial complex on $[n]$,
then it is known that $\alpha(u) \in K[x_1,\dots,x_n]$ for all $u \in G(\gin(I_\Gamma))$.
(See \cite[Lemma 8.15]{H}.)

\begin{theorem} \label{suff2}
Let $a$ and $b$ be positive integers.
One has $|\gins(I(K_{a,b}))|=1$.
\end{theorem}

\begin{proof}
Since $\overline {K_{a,b}} =K_a \cup K_b$ and $K_a \cup K_b$ is chordal,
it follows from Lemma \ref{froes} that $I_{\F (K_a \cup K_b)}=I(K_{a,b})$ has a linear resolution.
Then, by Lemma \ref{linear},
we have
\begin{eqnarray}
\beta_{ij}(\gin(I(K_{a,b})))=\beta_{ij}({\Gin(J_{\F(K_a\cup K_b)})}^*) \ \ \ \mbox{ for all } i,j. \label{SRM2}
\end{eqnarray}
On the other hand,
by Lemma \ref{finite},
we have 
$$\mathrm{min}_{\leq k} ( \Gin_{\mathrm{lex}}(J_{\F(K_a \cup K_b) }),2) = \mathrm{min}_{\leq k} ( \gin_{\mathrm{lex}}(I(K_{a,b})),2) \ \ \ \mbox{ for all } k.$$
Since $\min \{ t: x_t \mbox{ divides } u\} = \min \{ t: x_t \mbox{ divides } \alpha(u)\}$,
we have 
$$\mathrm{min}_{\leq k} ( \Gin_{\mathrm{lex}}(J_{\F(K_a \cup K_b) }),2) =
\mathrm{min}_{\leq k} ( \alpha(\gin_{\mathrm{lex}}(I(K_{a,b}))),2) \ \ \ \mbox{ for all }k.$$
Then, since $\alpha (\gin(I(K_{a,b}))$ is squarefree strongly stable,
Lemmas \ref{equivalent} says that
$\alpha(\gin_{\mathrm{lex}} (I(K_{a,b})))_2 ={ \Gin_{\mathrm{lex}}   (J_{\F(K_a \cup K_b)})}^*_2$.
Then the equations (\ref{SRM1}) say that
$$\beta_{ii+2} (\gin_{\mathrm{lex}}(I(K_{a,b})))= \beta_{ii+2} ( {\Gin_{\mathrm{lex}}(J_{\F(K_a \cup K_b)}  )}^*) \ \ \mbox{ for all } i.$$
We already proved $\Gin_{\mathrm{lex}}(J_{ \F(K_a \cup K_b) })= \Gin(J_{\F(K_a \cup K_b)})$ in Theorem \ref{extflag}.
Thus the above equation and (\ref{SRM2}) say that
$\beta_{ii+2}(\gin(I(K_{a,b})))= \beta_{ii+2} ( \gin_{\mathrm{lex}} (I({K_{a,b}})) )$ for all $i$.
Thus Lemma \ref{equivalent} and Corollary \ref{unique} say that
$\gin_{\mathrm{lex}}(I(K_{a,b}))_2=\gin(I(K_{a,b}))_2$ and $|\gins(I(K_{a,b}),2)|=1$.
Then, since Lemma \ref{froes} says that $\mathrm{reg}(I(K_{a,b}))=2$,
we have $|\gins(I(K_{a,b}))|=1$ by Lemma \ref{regularity} as required. 
\end{proof}

Next, we will prove the ``only if" part of Theorem \ref{main2}.

\begin{lemma} \label{4graphs}
Let $G$ be a graph.
If $|\gins(I(G))|=1$,
then $G$ contains none of the graphs $(a)$, $(b)$, $(c)$ and a cycle of length $3$
as an induced subgraph.
\end{lemma}

\begin{proof}
By Lemma \ref{equivalent}, we have $|\T(I(G),2)|=1$.
Then the claim for the graphs (a), (b) and (c) follows from Lemma \ref{syminduced} in the same way as Proposition \ref{threegraphs}.
We will show that if $G$ contains a cycle of length $3$ as an induced subgraph
then $|\gins(I(G))|>1$.
By Corollary \ref{unique} and Lemma \ref{syminduced},
what we must prove is that
if $I=(x_1x_2,x_1x_3,x_2x_3)$ then $|\T (I,2)|>1$.
Let 
$$I'= \init_{\mathrm{lex}}( \varphi_{2,3} (\init_{\mathrm{lex}} ( \varphi_{1,2} (I))))
\mbox{ and }I''=\init_{\mathrm{lex}}( \varphi_{1,3} (\init_{\mathrm{lex}} ( \varphi_{1,2} (I)))).$$
Then $I_2'$ and $I_2''$ are strongly stable.
Moreover, $I_2'$ is the $K$-vector space spanned by the monomials $x_1^2,x_1x_2,x_2^2$
and $I_2''$ is spanned by $x_1^2,x_1x_2,x_1x_3$.
Thus we have $|\T (I,2)|>1$.
\end{proof}

\begin{lemma} \label{symvarious}
Let $G$ be a graph on $[n]$ with $|\gins(I(G))|=1$.
Then
\begin{itemize}
\item[(i)] $G$ does not contain an induced cycle of length $\geq 5$;
\item[(ii)] If $G$ has more than  two proper connected components,
then $G$ has at most two edges;
\item[(iii)]
If $G$ is a connected bipartite graph, then $G$ is a complete bipartite graph.
\end{itemize}
\end{lemma}

\begin{proof}
The proofs for (i) and (iii) are the same as Lemma \ref{various}.
We will show (ii).
If $G$ has more than two proper connected components and have more than three edges,
since $G$ does not have an induced cycle of length $3$,
the graph $G$ must contain the graphs (b) or (c) as an induced subgraph.
Thus $G$ has at most two edges.
\end{proof} 

\begin{exam} \label{cocoa}
If $G$ has two proper connected components and has exactly two edges,
then we may assume that $I(G)=(x_1x_2,x_3x_4)$.
Then \cocoa's computation says that
$$\gin_{\mathrm{lex}}(I(G))=(x_1^2,x_1x_2,x_1x_3,x_2^4)$$
and
$$\gin(I(G))=(x_1^2,x_1x_2,x_2^3).$$
\cocoa\ computes the initial ideal $\init_\sigma(\varphi(I))$ for a random matrix $\varphi$.
Thus we can not guarantee that the above computations are true.
However, they are transformed strongly stable ideals of $I(G)$.
Then, by using the fact that the homogeneous component of degree $3$ of $(x_1^2,x_1x_2,x_1x_3,x_2^4)$
is spanned by all monomials $u$ of degree $3$ with $u \geq_{\mathrm{lex}} x_1x_3^2$,
Lemma \ref{conca} implies that it is in fact the same as the homogeneous component of degree $3$ of $\gin_{\mathrm{lex}}(I(G))$.
In particular, we have $x_2^3 \not \in \gin_{\mathrm{lex}}(I)$. 
On the other hand, by using the fact that $(x_1^2,x_1x_2,x_2^3)$ contains all monomials $u$
of degree $3$ with $u \geq_{\mathrm{rev}} x_2^3$,
Lemma \ref{conca} says that $x_2^3 \in \gin(I(G))$.
Thus we have $|\gins(I(G))|>1$.
\end{exam}

Then,
by the above example together with Lemma \ref{symvarious} (ii), we have
$|\gins(I(G))|>1$ if $G$ is a graph which has more than two proper connected components.

\begin{theorem} \label{suff}
Let $G$ be a graph on $[n]$.
If $|\gins(I(G))|=1$, then $G$ is a semi-complete bipartite graph.
\end{theorem}

\begin{proof}
By Lemma \ref{symvarious} (ii) together with Example \ref{cocoa},
it follows that $G$ has one proper connected component.
Also, by Lemmas \ref{4graphs} and \ref{symvarious} (ii),
every induced cycle of $G$ has length $4$.
Thus it follows from \cite[Proposition 1.6.1]{D} that
$G$ is a bipartite graph.
Then Lemma \ref{symvarious} (iii) says that $G$ is a complete bipartite graph.
\end{proof}


\begin{thebibliography}{1}

\bibitem{AHH} 
A. Aramova, J. Herzog and T. Hibi,
Ideals with stable Betti numbers,
\textit{Adv. Math. }\textbf{152} (2000) no. 1, 72--77.

\bibitem{SQlex} 
A. Aramova, J. Herzog and T. Hibi,
Squarefree lexsegment ideals,
\textit{Math. Z. } \textbf{122} (2000) no. 2, 353--378.

\bibitem{BNT}
E. Babson, I. Novik and R. Thomas,
Reverse lexicographic and lexicographic shifting,
\textit{J. Algebraic Combin.}, to appear.


\bibitem{coco}
CoCoA Team. 
\textit{CoCoA: a system for doing Computations in Commutative Algebra},
Available at {\tt http://cocoa.dima.unige.it.}


\bibitem{C}
A. Conca,
Reduction numbers and initial ideals,
\textit{Proc. Amer. Math. Soc. }\textbf{131} (2003) no. 4, 1015--1020.

\bibitem{C2}
A. Conca,
Koszul homology and extremal property of Gin and Lex,
\textit{Trans.  Amer. Math. Soc. }\textbf{256} (2004) no. 7, 2945--2961.




\bibitem{D}
R. Diestel,
``Graph theory, second edition",
Graduate text in Mathematics, \textbf{173},
Spriger--Verlag, New York, 2000.


\bibitem{E} D. Eisenbud,
``Commutative algebra, with a view toward algebraic geometry",
Graduate texts in Mathematics, \textbf{150},
Spriger--Verlag, New York, 1995.






\bibitem{F} 
R. Fr\"oberg, On Stanley--Reisner ring,
\textit{in} ``Topics in algebra",
Banach Center Publications, \textbf{26} Part 2, (1990), 57--70.


\bibitem{H} J. Herzog, Generic initial ideals and graded Betti numbers,
\textit{in} ``Computational Commutative Algebra and Combinatorics''
(T. Hibi, Ed.), Advanced Studies in Pure Math., Volume 33, (2002),
pp. 75--120.



\bibitem{Kh} G. Kalai,
Hyperconnectivity of graphs,
\textit{Graphs combin.}, \textbf{1}, (1985),  65--79.


\bibitem{K}
G. Kalai, Algebraic shifting,
\textit{in} ``Computational Commutative Algebra and Combinatorics''
(T. Hibi, Ed.), Advanced Studies in Pure Math., Volume 33, (2002),
pp. 121--163.


\bibitem{Mf} S. Murai,
Algebraic shifting of finite graphs,
arXiv:math.CO/0505010, \textit{comm. in Alg.}, to appear.

\bibitem{Mj} S. Murai,
Generic initial ideals and exterior algebraic shifting of join  of simplicial complexes,
arXiv:math.CO/0506298,
\textit{Arkiv f\"or matematik}, to appear.

\bibitem{HM1} S. Murai and T. Hibi,
The behavior of graded Betti numbers via algebraic shifting and combinatorial shifting,
arXiv:math.AC/0503685, (2005),  preprint.








\bibitem{N} E. Nevo,
Algebraic shifting and basic constructions on simplicial complexes,
\textit{J. Algebraic Combin.}, \textbf{22} (2005), 411--433.


\end{thebibliography}
\end{document}